\pdfoutput=1
\RequirePackage{ifpdf}
\ifpdf 
\documentclass[pdftex]{sigma}
\else
\documentclass{sigma}
\fi

\numberwithin{equation}{section}

\newtheorem{Theorem}{Theorem}[section]
\newtheorem*{Theorem*}{Theorem}

\newtheorem{Lemma}[Theorem]{Lemma}
\newtheorem{prop}[Theorem]{Proposition}

\newtheorem{assumption}[Theorem]{Assumption}

\theoremstyle{definition}
\newtheorem{Definition}[Theorem]{Definition}

\newtheorem{Remark}[Theorem]{Remark}

\usepackage{tikz-cd}

\newcommand{\opP}{\operatorname{P}}

\newcommand{\fraks}{\mathfrak{s}}

\newcommand{\RR}{\mathbb{R}}
\newcommand{\R}{\mathbb{R}}
\newcommand{\mcK}{\mathcal{K}}
\newcommand{\dd}{\mathrm{d}}
\newcommand{\cc}{\mathrm{C}}
\newcommand{\mcB}{\mathcal{B}}

\newcommand{\mfs}{\mathfrak{s}}
\newcommand{\bbb}{\mathrm{B}}

\newcommand{\eps}{\varepsilon}

\newcommand{\E}{\mathbb{E}}

\newcommand{\bfK}{\pmb{\mathcal{K}}}

\def\enlarge{edge[draw=none] (0,-0.1) edge[draw=none] (0,0.1)}

\colorlet{symbols}{blue!90!black}
\colorlet{testcolor}{green!60!black}

\colorlet{testcolor}{green!60!black}

\usetikzlibrary{calc}
\usetikzlibrary{shapes.misc}
\usetikzlibrary{shapes.symbols}
\usetikzlibrary{shapes.geometric}
\usetikzlibrary{snakes}
\usetikzlibrary{decorations}
\usetikzlibrary{decorations.markings}

\tikzset{
	root/.style={circle,fill=testcolor,inner sep=0pt, minimum size=2mm},
	broot/.style={circle,fill=gray,inner sep=0pt, minimum size=2mm},
	dot/.style={circle,fill=black,inner sep=0pt, minimum size=1mm},
		reddot/.style={circle,fill=red,inner sep=0pt, minimum size=1mm},
			bluedot/.style={circle,fill=blue,inner sep=0pt, minimum size=1mm},
	eps/.style={circle,fill=white,draw=symbols,inner sep=0pt,minimum size=0.8mm},
	int/.style={circle,fill=black,draw=black,inner sep=0pt,minimum size=0.7mm},
	var/.style={circle,fill=black!10,draw=black,inner sep=0pt, minimum size=2mm},
	dotred/.style={circle,fill=black!50,inner sep=0pt, minimum size=2mm},
	generic/.style={semithick,shorten >=1pt,shorten <=1pt},
	dist/.style={ultra thick,draw=testcolor,shorten >=1pt,shorten <=1pt},
	testfcn/.style={ultra thick,testcolor,shorten >=1pt,shorten <=1pt,<-},
	testfcnx/.style={ultra thick,testcolor,shorten >=1pt,shorten <=1pt,<-,
		postaction={decorate,decoration={markings,mark=at position 0.6 with {\drawx}}}},
	keps/.style={semithick,shorten >=1pt,shorten <=1pt,densely dashed,->},
	kprimex/.style={semithick,shorten >=1pt,shorten <=1pt,densely dashed,->,
		postaction={decorate,decoration={markings,mark=at position 0.4 with {\drawx}}}},
	kernel/.style={semithick,shorten >=1pt,shorten <=1pt,->},
	multx/.style={shorten >=1pt,shorten <=1pt,
		postaction={decorate,decoration={markings,mark=at position 0.5 with {\drawx}}}},
	kernelx/.style={semithick,shorten >=1pt,shorten <=1pt,->,
		postaction={decorate,decoration={markings,mark=at position 0.4 with {\drawx}}}},
	kernel1/.style={->,semithick,shorten >=1pt,shorten <=1pt,postaction={decorate,decoration={markings,mark=at position 0.45 with {\draw[-] (0,-0.1) -- (0,0.1);}}}},
	kernel2/.style={->,semithick,shorten >=1pt,shorten <=1pt,postaction={decorate,decoration={markings,mark=at position 0.45 with {\draw[-] (0.05,-0.1) -- (0.05,0.1);\draw[-] (-0.05,-0.1) -- (-0.05,0.1);}}}},
	kernelBig/.style={semithick,shorten >=1pt,shorten <=1pt,decorate, decoration={zigzag,amplitude=1.5pt,segment length = 3pt,pre length=2pt,post length=2pt}},
	rho/.style={dotted,semithick,shorten >=1pt,shorten <=1pt},
	renorm/.style={shape=circle,fill=white,inner sep=1pt},
	labl/.style={shape=rectangle,fill=white,inner sep=1pt},
	xi/.style={circle,fill=symbols!10,draw=symbols,inner sep=0pt,minimum size=1.2mm},
	xix/.style={crosscircle,fill=symbols!10,draw=symbols,inner sep=0pt,minimum size=1.2mm},
	xib/.style={circle,fill=symbols!10,draw=symbols,inner sep=0pt,minimum size=1.6mm},
	xibx/.style={crosscircle,fill=symbols!10,draw=symbols,inner sep=0pt,minimum size=1.6mm},
	not/.style={circle,fill=symbols,draw=symbols,inner sep=0pt,minimum size=0.5mm},
cumu2n/.style={inner sep=3pt},
cumu2/.style={draw=red!80,fill=red!40},
cumu2b/.style={draw=blue!80,fill=blue!40},
cumu2nv/.style={inner sep=3pt},
cumu2v/.style={draw=red!80,fill=white,very thick},
cumu3/.style={regular polygon, regular polygon sides=3,draw=red!80,rounded corners=3pt,fill=red!40,minimum size=5mm},
cumu4/.style={regular polygon, regular polygon sides=4,draw=red!80,rounded corners=3pt,fill=red!40,minimum size=7mm},
cumu5/.style={regular polygon, regular polygon sides=5,draw=red!80,rounded corners=3pt,fill=red!40,minimum size=7mm},
	>=stealth,
	not/.style={circle,fill=symbols,draw=symbols,inner sep=0pt,minimum size=0.5mm},
kernels2/.style={very thick,segment length=12pt},
	}
\makeatletter
\def\DeclareSymbol#1#2#3{%
	\expandafter\gdef\csname MH@symb@#1\endcsname{\tikzsetnextfilename{symbol#1}%
		\tikz[baseline=#2,scale=0.15,draw=symbols,line join=round]{#3}}%
	\expandafter\gdef\csname MH@symb@#1s\endcsname{\scalebox{0.75}{\tikzsetnextfilename{symbol#1}%
			\tikz[baseline=#2,scale=0.15,draw=symbols,line join=round]{#3}}}%
	\expandafter\gdef\csname MH@symb@#1ss\endcsname{\scalebox{0.65}{\tikzsetnextfilename{symbol#1}%
			\tikz[baseline=#2,scale=0.15,draw=symbols,line join=round]{#3}}}%
}
\def\<#1>{\ifthenelse{\boolean{mmode}}{\mathchoice{\csname MH@symb@#1\endcsname}{\csname MH@symb@#1\endcsname}{\csname MH@symb@#1s\endcsname}{\csname MH@symb@#1ss\endcsname}}{\csname MH@symb@#1\endcsname}}
\makeatother

\tikzset{
	eps/.style={circle,fill=white,draw=symbols,inner sep=0pt,minimum size=0.8mm},
	}
\makeatletter
\def\DeclareSymbol#1#2#3{\expandafter\gdef\csname MH@symb@#1\endcsname{\tikz[baseline=#2,scale=0.15,draw=symbols]{#3}}}
\def\<#1>{\csname MH@symb@#1\endcsname}
\makeatother

\DeclareSymbol{1}{0}{\draw[white] (-.5,0) -- (.5,0); \draw (0,0) -- (0,1.5) node[eps] {};}
\DeclareSymbol{2}{0}{\draw (-0.5,1.5) node[eps] {} -- (0,0) -- (0.5,1.5) node[eps] {};}
\DeclareSymbol{3}{0}{\draw (0,0) -- (0,1.5) node[eps] {}; \draw (-.7,1.3) node[eps] {} -- (0,0) -- (.7,1.3) node[eps] {};}

\DeclareSymbol{b1}{0}{
\draw[white] (-.5,0) -- (.5,0); \draw[kernels2] (0,0) -- (0,1.5);
 \draw (0,1.5) node[eps] {};
}

\DeclareSymbol{b2}{0}{
\draw[kernels2] (0,0) -- (-0.5,1.4);
\draw[kernels2] (0,0) -- (+0.5,1.4);
 \draw (-0.5,1.4) node[eps] {};
 \draw (0.5,1.4) node[eps] {};
 }

\makeatletter

\DeclareRobustCommand{\TitleEquation}[2]{\texorpdfstring{\StrLeft{\f@series}{1}[\@firstchar]$\if%
		b\@firstchar\boldsymbol{#1}\else#1\fi$}{#2}}

\makeatother

\def\CCE{\mathbb{E}}
\def\CCV{\mathbb{V}}
\def\CCG{\mathbb{G}}

\DeclareSymbol{Xi2}{0}{\draw[white] (-.5,0) -- (.5,0); \draw (0,0) -- (0,1.5) node[eps] {};
\draw (0,0) node[eps] {};}

\begin{document}

\allowdisplaybreaks

\newcommand{\arXivNumber}{2507.23737}

\renewcommand{\thefootnote}{}

\renewcommand{\PaperNumber}{021}

\FirstPageHeading

\ShortArticleName{Renormalisation of Singular SPDEs with Correlated Coefficients}

\ArticleName{Renormalisation of Singular SPDEs\\ with Correlated Coefficients\footnote{This paper is a~contribution to the Special Issue on Asymptotics, Randomness and Noncommutativity. The~full collection is available at \href{https://sigma-journal.com/noncommutativity.html}{https://sigma-journal.com/noncommutativity.html}}}

\Author{Nicolas CLOZEAU~$^{\rm a}$ and Harprit SINGH~$^{\rm b}$}

\AuthorNameForHeading{N.~Clozeau and H.~Singh}

\Address{$^{\rm a)}$~Universit\'e de Toulon, Av. de l'Universit\'e, 83130 La Garde, France}
\EmailD{\mail{clozeau@univ-tln.fr}}
\URLaddressD{\url{https://sites.google.com/view/nicolasclozeauphd}}

\Address{$^{\rm b)}$~SISSA, Via Bonomea 265, 34136 Trieste TS, Italy}
\EmailD{\mail{harprit.singh@bluewin.ch}}
\URLaddressD{\url{https://sites.google.com/view/harprit-singh/home}}

\ArticleDates{Received September 01, 2025, in final form February 20, 2026; Published online March 07, 2026}

\Abstract{We show local well-posedness of the g-PAM and the \smash{$\phi^{K+1}_2$}-equation for $K\geq 1$ on the two-dimensional torus when the coefficient field is random and correlated to the driving noise. In the setting considered here, even when the model in the sense of Hairer (2014) is stationary, naive use of renormalisation constants in general leads to variance blow-up. Instead, we prove convergence of renormalised models choosing random renormalisation functions analogous to the deterministic variable coefficient setting. The main technical contribution are stochastic estimates on the model in this correlated setting which are obtained by a combination of heat kernel asymptotics, Gaussian integration by parts formulae and Hairer--Quastel type bounds.}


\Classification{60H17; 60L30}

\renewcommand{\thefootnote}{\arabic{footnote}}
\setcounter{footnote}{0}

\section{Introduction}

 For a random, uniformly elliptic coefficient field $a=\{a_{ij}\}_{1\leq i,j\leq 2}$,
 we consider the generalised parabolic Anderson model (g-PAM) equation, formally given by
\begin{gather}\label{eq:PAM intro}
\partial_{t}u- \sum_{i,j=1}^2 a_{ij}(x) \partial_{i }\partial_{j} u = \sum_{i,j=1}^2 f_{ij}(u) \partial_i u\, \partial_j u + g(u)\xi\qquad\text{on $(0,\infty)\times \mathbb{T}^2$},
\end{gather}
where $\xi$ denotes a spatial white noise on $\mathbb{T}^2$; and the \smash{$\phi^{K+1}_2$}-equation (for any $K\geq 1$) which formally reads
\begin{gather}\label{eq:phi^4}
\partial_{t}u- \sum_{i,j=1}^2 a_{ij}(t,x) \partial_{i }\partial_{j} u = -u^{K} + \xi \qquad\text{on $(0,\infty)\times \mathbb{T}^2$},
\end{gather}
where $\xi$ denotes a space-time white noise on $(0,\infty)\times \mathbb{T}^2$.
While the local well-posedness and renormalisation of these equations are well understood for deterministic, sufficiently regular coefficient fields,
 we study here a setting where $a$ is correlated with the driving noise on the right-hand side, see Assumption~\ref{AssumCoef}.

The g-PAM equation with constant coefficients was first solved in the pioneering\footnote{These articles initiated extensive progress in the understanding of subcritical singular SPDEs, too vast to be accounted for here. Instead we refer to the textbook~\cite{frizhairerbook} and the references therein.} works \cite{GIP15, Hai14} and forms an important example of a singular SPDE.
On the one hand, it is a generalisation of the (continuum) linear PAM, cf.\ \cite{gartner2005parabolic, HL15, HL18, hsu24, KPZ22, MP19, Mat22}, a model of a population/branching process describing the expected density of diffusing particles in the presence of a random environment of sources and sinks given by the noise. The variable
coefficient field is then another environment governing the diffusivity of these particles, which would make the assumption that the two environments are independent rather artificial. On the other hand, this equation has been a popular test case in the study of singular SPDEs in geometries beyond the well studied setting of constant coefficient equations in flat spaces, cf.\ \cite{BB16, BAUD25, BAUD, DDD19, HS23m, MS23, Mou22}.

The $\phi_2^{4}$-equation with constant coefficients was first solved in~\cite{DD03}, see also~\cite{TW18}, and is a~Langevin-dynamic for the associated (Gibbs)-measure, see~\cite{GlimmJaffe, Simon}. It can also be interpreted as a toy model for the dynamics of a ferromagnet close to the critical temperature, cf.~\cite{mourrat2017convergence}, in~which case the coefficient field can be seen as inhomogeneities/random bonds in the material, which we allow to be correlated to the external forcing.

While this article is certainly not the first to consider SPDEs with random, correlated coefficients, see for instance \cite{Car99, Ferrante2006spdes}, we are not aware of any previous works considering the type of coefficients studied here for singular equations.
The main novelty for singular SPDEs with correlated coefficients, compared to the well studied deterministic coefficient setting, is in establishing the required stochastic estimates as will be discussed below.

\subsection{Solution theory based on Rough Analysis}
Building on the ideas originating in rough path theory \cite{gub04, lyo98}, there are by now several approaches to tackle singular SPDEs such as regularity structures
\cite{Hai14}, paracontrolled calculus \cite{GIP15} or renormalisation group approaches \cite{duc25, Kup16}. Roughly speaking, in all of these approaches there are two steps to perform, which we here perform within the framework of regularity structures.%
\begin{enumerate}\itemsep=0pt
\item Consider the solution to the linear equation and constructs a finite number of (appropriately renormalised) functions thereof.
This data, roughly speaking, provides the so-called~\textit{model}.
\item Given a model, one solves an enhanced PDE.
\end{enumerate}
For the second step, we will essentially be able to refer to \cite{Hai14, Sin25, TW18}, which allows us to mainly focus on the novelties appearing due to the presence of correlations and outsource aspects of the problem that are well understood.

Next, we make precise the assumed structure and regularity of the coefficient field, where the
 regularity assumption is used through properties of the parabolic Green's function, see also Remark~\ref{rem:regularity} below.

\begin{assumption}[structural assumption on the coefficient field]\label{AssumCoef}
For $\lambda>0$, let $A\colon \mathbb{R}\to \mathbb{R}^{2\times 2}$ take values in the set of $\lambda$-uniformly elliptic\footnote{This means that
$A(s)\zeta\cdot \zeta\geq \lambda\vert \zeta\vert^2$ uniformly over $s\in\mathbb{R}$ and $\zeta\in \mathbb{R}^{2\times 2}$.}
 matrices satisfying
\[\sup_{s\in \mathbb{R}}\biggl\vert\frac{\dd^k}{\dd s^k}A(s)\biggr\vert<\infty\qquad\text{for any $k\geq 0$.}\]
Depending on the equation, we make the following structural assumption:
\begin{itemize}\itemsep=0pt
\item
For g-PAM, we consider a spatial white noise $\xi\in \mathcal{D}'\bigl(\mathbb{T}^2\bigr)$
and for $\mu, \sigma\in C^2\bigl(\mathbb{T}^2\bigr)$ set
$ h:=\sigma*\xi+\mu$.

\item For the \smash{$\phi^{K+1}_2$}-equation, we consider a space-time white noise $\xi\in \mathcal{D}'\bigl(\mathbb{R}\times \mathbb{T}^2\bigr)$
and for\footnote{The compact support assumption on $\sigma$ is only non-empty in the time direction and could be replaced by an~appropriate integrability.} \smash{$\mu\in C^3 \bigl(\mathbb{R}\times \mathbb{T}^2\bigr)$}, \smash{$\sigma \in C_c^3 \bigl(\mathbb{R}\times \mathbb{T}^2\bigr)$} set
$ h:=\sigma*\xi+\mu $
 \big(where $*$ denotes the space-time convolution on $\mathbb{R}\times \mathbb{T}^2$\big).
\end{itemize}
In both cases, we build our coefficient field as
$a=A(h)$ which, in the former case, is a coefficient field $a\colon \mathbb{T}^{2}\to \mathbb{R}^{2\times 2}$
and in the latter case $a\colon \mathbb{R}\times \mathbb{T}^2 \to \mathbb{R}^{2\times 2}$.
\end{assumption}

\begin{Remark}In both examples, the coefficient field has a direct interpretation as a random environment/heterogeneities in materials, see \cite{Torquato_book}. The study of \eqref{eq:PAM intro} and \eqref{eq:phi^4} is therefore directly motivated by the relation of these SPDEs to statistical mechanics.
 The specific choice to take coefficients of form $A(h)$ is inspired by examples in \cite{AKMbook, Homogenisation_GNO}.
For the parabolic Anderson model, even the case $\sigma=1$, $\mu=0$ is novel, and requires the use of random renormalisation constants (as discussed below).
\end{Remark}

Recall, that in order to solve the linear parabolic Anderson model, i.e., \eqref{eq:PAM intro} with $f_{i,j}=0$ and~${g(u)=u}$, constructing the model essentially\footnote{Note that in regularity structures one does not quite work with the solution of the equation directly, but the same mechanism also leads to variance blow up for the model.
} amounts to defining an appropriately renormalised
product between the solution $u\colon (0,\infty)\times \mathbb{T}^2 \to \mathbb{R}$ to the equation
\begin{equation}\label{SolutionStochasticHeatEq}
\partial_t u - \sum_{i,j=1}^2 a_{ij}(x) \partial_i \partial_j u= \xi, \qquad u(0)=0,
\end{equation}
and the spatial white noise $\xi$ on $\mathbb{T}^2$ (interpreted as a function of space-time). In order to formulate results involving spatial white noise, we shall work with mollifications built from a~function~${\rho\in \cc^\infty_c(\bbb_1)}$ such that, for any $x\in \mathbb{T}^2$, $\rho(-x)=\rho(x)$ and $\int_{\mathbb{R}^2} \rho(x)\, \dd x=1$ by setting, for $\delta\in (0,1]$
\begin{gather}\label{eq:regularised noise}
\xi_\delta = \xi*\rho^\delta \qquad\text{for}\ \rho^\delta:=\delta^{-2}\rho\bigl(\tfrac{\cdot}{\delta}\bigr).
\end{gather}
Considering the case when $\mu=0$ in Assumption~\ref{AssumCoef}, both random fields $\xi$ and $u$ are spatially stationary,
 potentially suggesting renormalisation by deterministic constants as is common in the singular SPDE literature.
 The next proposition, the proof of which can be found in Section~\ref{sec:proof of}, in particular shows that this is not possible as it
leads to either mean or variance blow up.

\begin{prop}[variance blow-up]\label{prop:varblowup}
In the setting of Assumption~{\rm \ref{AssumCoef}} for the g-PAM equation, assume furthermore that
$\sigma \neq 0$ and $\det (A)$ is not constant.
Let $u_\delta$ be the solution to~\eqref{SolutionStochasticHeatEq} with the white noise $\xi$ replaced by $\xi_\delta$ defined in~\eqref{eq:regularised noise}. Then, for any deterministic sequence of functions $\{c_\delta\}_{\delta>0}\subset C\bigl(\mathbb{T}^2\bigr)$ at least one of the following two assertions holds.
\begin{itemize}\itemsep=0pt
\item[$(i)$] There exists \smash{$\phi\in \cc^\infty\bigl(\mathbb{T}^2\bigr)$} such that
$\limsup_{\delta\to 0} | \mathbb{E}[(u_\delta(1,\cdot)\xi_\delta-c_\delta,\phi)]|= \infty$.
\item[$(ii)$]There exists \smash{$\phi\in \cc^\infty\bigl(\mathbb{T}^2\bigr)$} such that
$\limsup_{\delta\to 0} \mathbb{E}\bigl[(u_\delta(1,\cdot)\xi_\delta-c_\delta,\phi)^2\bigr]= \infty$.
\end{itemize}
\end{prop}

 This type of variance blow up, though due to a distinct mechanism, is well known in rough path theory, see \cite{Unt08} as well as the more recent works on SPDEs \cite{GT25, Hai24}.
 \begin{Remark}
 As will become clear in the proof, the phenomenon captured in Proposition~\ref{prop:varblowup} is generic and does also hold for the other stochastic object appearing for the g-PAM equation and those for the \smash{$\phi^{K+1}_2$}-equation.
 \end{Remark}

\subsection{Main results}

Our first main theorem, Theorem~\ref{thm:g-pam}, establishes local well posedness for the $g$-PAM equation~\eqref{eq:PAM intro}. We first define the involved renormalisation functions, which in contrast to the most studied case of constant coefficients singular SPDEs are genuine functions and not constants. To motivate our choice of functions, we recall that the Green's function $\Gamma(t,x,s,y)$ of the heat operator $\partial_t-\sum a_{ij}(x) \partial_i\partial_j$, by Levy's construction of fundamental solutions, has an explicit series representation, cf.~\cite{Fri08}. It was noted in~\cite{Sin25} that one can renormalise the g-PAM equation by neglecting all but the leading order term
\begin{gather}\label{eq:frozen kernel}
Z^{x}(t,y)= \frac{1}{4\pi t \sqrt{\det(a(x))}} \exp\biggl(- \frac{ y \cdot a^{-1}(x) y}{4t}\biggr) \qquad \text{for $x,y\in \mathbb{R}^2$, $t\in (0,\infty)$}
\end{gather}
in that representation.
Analogously to \cite[Section~3.1.1]{Sin25}, we define for $x\in \mathbb{R}^2$ the kernels on~$\mathbb{R}^2\setminus \{ 0\}$
\begin{gather}\label{eq:space dep frozen}
G^{x}(y)= \int_{0}^1 Z^{x}(t,y)\, \dd t, \qquad
G_i^{x}(y)= \sum_{j=1}^{2} a^{ij}(x) y_j\int_{0}^1 \frac{1}{2 t}Z^{x}(t,y)\, \dd t,
\end{gather}
where
$a^{ij}:=\bigl(a^{-1}\bigr)_{ij}$. Finally, given a mollifier $\rho\colon \mathbb{R}^2 \to \R$ as for \eqref{eq:regularised noise}
the random renormalisation functions are chosen as
\begin{gather}
c_{\delta}^{\<Xi2>}(x) = \int_{\mathbb{R}^2} G^{x}(y) \rho^\delta *\rho^\delta (y)\, \dd y,\nonumber \\
c_{\delta}^{\<b2>_{ij}}(x) = \int_{\mathbb{R}^2\times \mathbb{R}^2} G_i^{x}(y)G_j^{x}(y') \rho^\delta *\rho^\delta(y-y')\, \dd y \dd y'.\label{eq:renormalisation function pam}
\end{gather}

\begin{Remark}
Observe that this choice of counterterms is such that in the case $\sigma=0$ the
solutions agree exactly with the ones constructed in \cite{Sin25}. For the motivation of that choice of counterterms, see the discussion in Section~1.1 therein.
One can directly compute that
\[
c_\delta^{\<Xi2>}, \qquad c^{\<b2>}_\delta=\bigl\{c^{\<b2>_{ij}}_\delta\bigr\}_{ij}
\]
in the theorems below can be written as
\[c_\delta^{\<Xi2>}(x) = \frac{|\log (\delta)|}{2\pi \sqrt{\det(a(x))}}+ \beta_{\delta}^{\<Xi2>} (x), \qquad
 c^{\<b2>}_\delta = \frac{a^{-1}(x)|\log (\delta)|}{4\pi \sqrt{\det(a(x))}} +\beta_{\delta}^{\<b2>} (x), \]
where \smash{$\beta_{\delta}^{\<Xi2>} (\rho)$}, \smash{$\beta_{\delta}^{\<b2>} (\rho)$} are some $\rho$-dependent functions and converge as $\delta \to 0$, cf.\ \cite[Remark~1.18]{HSper}.
\end{Remark}

\begin{Remark}
Heuristically, the fact that the renormalisation functions for the considered equations should depend on the coefficient field in a local way can also be argued for similarly to \cite[Section~1.5]{broux}.
\end{Remark}

\begin{Theorem}\label{thm:g-pam}
Let $a$ and $\xi$ be as in Assumption~{\rm \ref{AssumCoef}}, and let $\xi_\delta$ be its regularisation for a~mollifier~$\rho$ as in \eqref{eq:regularised noise}. Let $f_{ij}$, $g\in C^{\infty}(\RR)$ and $u_0\in C^\alpha\bigl(\mathbb{T}^2\bigr)$ for $\alpha>0$.
There exist a~random~${T>0}$ and, for each $\delta\in (0,1]$, a random process $u_\delta\in C\bigl([0,T]\times \mathbb{T}^2\bigr)$ satisfying $u_\delta(0)=u_0$ and $($in the mild sense$)$
\begin{gather}
 \partial_t u_\delta- \sum_{i,j=1}^d a_{ij}(x) \partial_i \partial_j u_\delta = \sum_{i,j=1}^2 f_{ij}(u_\delta) \bigl( \partial_i u_\delta\,\partial_j u_\delta - c^{\<b2>_{ij}}_\delta (x) g^2(u_\delta)\bigr) \nonumber\\
 \hphantom{\partial_t u_\delta- \sum_{i,j=1}^d a_{ij}(x) \partial_i \partial_j u_\delta =}{}
 + g(u_\delta)\bigl( \xi_\delta -c_{\delta}^{\<Xi2>}(x) g'(u_\delta)\bigr)\label{eq:PAM}
\end{gather}
on $(0,\infty)\times \mathbb{T}^2$.
Furthermore, there exists $u\in C\bigl([0,T]\times \mathbb{T}^2\bigr)$ independent of the choice of $\rho$ such that $u_\delta \to u$ uniformly on $[0,T]\times \mathbb{T}^2$ as $\delta\to 0$ in probability.
\end{Theorem}
\begin{Remark}
We expect that the solutions constructed here agree with the solutions one would obtain by regularising the noise `covariantly' as in \cite[Section~3.1]{Sin25} and subtracting
an exact multiple of \smash{$\det(a(t,x))^{-1/2}$}, resp.\ \smash{$\det(a(t,x))^{-1/2} a^{-1}(t,x)$}, see Remark~4.8 therein.
The reason we do not pursue this here is that it would require a modification of the Hairer--Quastel criterion, the proof of which would be quite lengthy to formulate.
\end{Remark}

Next, we present the well-posedness theorem for the $\phi^{K+1}_2$-equation \eqref{eq:phi^4}. This time the noise~$\xi$ denotes a space-time white noise, and we shall typically
write $z=(s,x)\in \mathbb{R}\times\mathbb{T}^2$. We~shall work here with space time mollifications: for $\rho\in \cc^\infty_c(B_1)$ satisfying $\rho(-z)=\rho(z)$ and $\int \rho=1$, we set for $\delta\in (0,1)$
\begin{gather}\label{eq:space time mollified noise}
\xi_{\delta}= \xi*\rho^\delta \qquad\text{for}\ \rho^{\delta}(s,x)= \delta^{-4}\rho\bigl(\tfrac{s}{\delta^2},\tfrac{x}{\delta}\bigr),
\end{gather}
where this time the convolution is a space time convolution. Similarly to
\eqref{eq:frozen kernel}, set
\[Z^{z}(t,y)= \frac{1}{4\pi t \sqrt{\det(a(z))}} \exp\biggl( - \frac{ y \cdot a^{-1}(z) y }{4t} \biggr)\]
and motivated as before, define
\begin{gather}\label{eq:counterterm phi}
c_{\delta}^{\<2>}(t,x)= \int_{([0,1]\times \mathbb{R}^2)^{2}} Z^{(t,x)}(z)Z^{(t,x)}(z') \rho^\delta * \rho^\delta (z-z') \dd z \dd z',
\end{gather}
cf.~\cite[Section~3.3]{Sin25}. We denote by $H_N\colon \mathbb{R}\times\mathbb{R}\to \mathbb{R}$ the $N$-th generalised Hermite polynomial which is defined by the recursive relations
\begin{gather}\label{eq:Hermit}
H_N(X,C)=X H_{N-1}(X,C)- (N-1) C H_{N-2}(X,C),
\end{gather}
where $H_0=1$ and $H_1(X,C)=X$.

\begin{Theorem}\label{thm:phi4_2}
Let $a$ and $\xi$ be as in Assumption {\rm \ref{AssumCoef}}, let $\xi_\delta$
be its regularisation for a~mollifier~$\rho$ as in
 \eqref{eq:space time mollified noise}, and let \smash{$\bigl\{c_{\delta}^{\<2>}\bigr\}_{\delta>0}$} be as in \eqref{eq:counterterm phi}.
Let
$K\geq 0$ and
 \smash{$u_0\in C^{\alpha}\bigl(\mathbb{T}^d\bigr)$} for \smash{$\alpha>-\frac{1}{10}$}.
 There exist a random $T>0$ and for each $\delta\in (0,1]$ a random process \smash{$u_\delta\in C\bigl((0,T)\times \mathbb{T}^2\bigr)$} satisfying $u_\delta(0)=u_0$ and $($in the mild sense$)$
\begin{gather}\label{phi4}
\partial_t u_\delta -\sum_{i,j=1}^{2} a_{ij}(t,x)\partial_{ij} u_\delta= -H_{K}\bigl(u_\delta, c_{\delta}^{\<2>}\bigr)+\xi_\delta, \qquad
u_\delta(0)=u_0.
\end{gather}
Furthermore, there exists $u\in \mathcal{D}'\bigl((0,T)\times \mathbb{T}^2\bigr)$ independent of the choice of $\rho$ such that $u_\delta \rightarrow u$ as $\delta\rightarrow 0$ in probability.
\end{Theorem}

\begin{Remark}
Let us emphasise that, since the sequence \smash{$\bigl\{c_{\delta}^{\<2>}\bigr\}_{\delta>0}$} depends on the realisation of the underlying Gaussian noise $\xi$, the Wick polynomial \smash{$H_{K}\bigl(u_\delta, c^{\<2>}_\delta\bigr)$} does not coincide with the~$K$-th Wick power of the solution to the linear stochastic heat equation.
\end{Remark}
\begin{Remark}
Note that for $K$ odd, the random time $T$ in Theorem~\ref{thm:phi4_2} could be made deterministic and arbitrarily large following arguments from~\cite{MWplane,TW18}. Since there would be little novelty in this step, we do not perform it here. Similarly, the distribution space topology in which $u_\delta$ converges as $\delta\to 0$ could be made sharper.
\end{Remark}

As already mentioned, much of the work in proving Theorems~\ref{thm:g-pam} and~\ref{thm:phi4_2} is in obtaining almost sure bounds and convergence of the model.
These bounds for \eqref{eq:PAM intro} and \eqref{phi4} in the case of deterministic coefficients, using equivalence of moments for random variables in a finite Wiener chaos, cf.\ \cite[Lemma~10.5]{Hai14}, reduce to elementary computations. Here, due to the randomness in the coefficient field the model does not take values in a finite Wiener chaos, which requires us to use slightly more sophisticated tools. In order to apply the Kolmogorov criterion to pass from moment bounds to a.s.\ estimates we bound arbitrary moments of the model as for example also done in~\cite{CS17, HShen} in a different situation. Here, we first obtain an explicit integral representation of the $q$-th moment by the classical Isserlis theorem, see Appendix~\ref{App:Gaussian int}, and then bound each term in this representation using the Hairer--Quastel criterion \cite{HQ18} and a slight variant thereof, see Appendix~\ref{App:HairerQuastel}.

{\bf Scope for generalisation and outlook.}
The equations considered here are arguably the `simplest' singular SPDEs, but already for slightly more singular equations such as the $\phi^4_3$\nobreakdash-equation establishing the analogous stochastic estimates by the same tools would be extremely tedious. Thus, we believe it would be desirable to develop more systematic tools such as \cite{BM23,CH16, HS23, LOTT24} in the constant coefficient setting and \cite{BSS25} in the deterministic variable coefficient setting.

It would also be interesting to weaken the rather strong Assumption~\ref{AssumCoef}, on the one hand to reduce the regularity assumption on the matrix $A$ and on the other hand to more general correlation structures.

Finally, we view this article as a step-$0$ towards the study of stochastic homogenisation of singular SPDEs as it provides a reference solution theory to make precise homogenisation statements (such as in the periodic setting the relationship of \cite{Sin25} to \cite{HSper} and to the forthcoming~\cite{HSper2}, see also \cite{CFX24, CX23} for further results on periodic homogenisation of singular SPDEs).

\begin{Remark}\label{rem:regularity}
Note that the required regularity of the functions $\mu$, $\sigma$ in Assumption~\ref{AssumCoef} enters due to the used regularity assumption on the coefficient field $a$, which in turn enters only through estimates on the fundamental solution and bounds on the summands of its series expansion.
For the parabolic Anderson model, the slightly stronger then expected assumption that ${a\in C^2}$ is used
 to work with the kernel spaces $\bfK^{2}_{L,R}$ for $L,R>1$ as defined in \cite[Defenition~2.12]{Sin25}. Replacing these spaces with their slight variants introduced in the forthcoming article \cite{HSper2} would allow without any modification in the argument presented here to take $a\in C^{1,\alpha}$ for $\alpha >0 $ in~Theorem~\ref{thm:g-pam}.

For the $\phi^{K+1}_2$-equation, the strong regularity assumption is used in order to treat the remainder term \eqref{eq:higher regularity used} in the proof of Lemma~\ref{lem:phi4 estimate} by classical Young multiplication.\footnote{Recall that the Young multiplication theorem states that whenever $\alpha+\beta>0$ pointwise multiplication between two smooth functions extends uniquely to a continuous map $C^\alpha\times C^\beta\to C^{(\alpha\wedge\beta)}$.}
We believe that this assumption could be weakened to $a\in C^{\alpha}$ for arbitrary $\alpha>0$, but this would require bounds on the Malliavin derivative of the fundamental solution $\Gamma$. As this technical point would lengthen the article notably, we decided to not pursue this direction in the present contribution.
Such bounds could be obtained by, for example, Malliavin differentiating the PDE to which $\Gamma$ is the fundamental solution. We refer to \cite{clozeau2025inductive}, where this idea is used to derive stochastic estimates.
\end{Remark}

\begin{Remark}
Continuing the discussion of Remark~\ref{rem:regularity},
let us observe that another place where the stronger regularity assumption at first sight seems to be used in the article is when applying the Hairer--Quastel criterion, Theorem~\ref{thm:HQ}, resp.\ the variant Theorem~\ref{thm:slight variant of HQ} thereof in~Lemma~\ref{lem:new_term}, resp.\ Lemma~\ref{lem:phi4 singular estimate}. But in both cases this is easily circumvented by using H\"older seminorms instead of the supremum norm on the derivative to bound increments in the multi-scale clustering argument for the respective singular integrals
in Lemmas~\ref{lem:new_term}, \ref{lem:dumbbell2},~\ref{lem:cherry} and~\ref{lem:phi4 singular estimate}.
\end{Remark}

\subsection{Notations}
We shall always work on an ambient probability space $(\Omega, \mathcal{F}, P)$ supporting the white noise. For~${q\in [1,\infty]}$, we shall denote by $L^{q}(\Omega)$ scalar valued random variables with finite $q$-th moment. For two random variables $X,Y$ we write ${\rm Cov}( X,Y)= \E[(X-\E[X]) (Y-\E[Y])]$.

We identify the two dimensional torus $\mathbb{T}^2$ with the quotient $\mathbb{R}^2/\mathbb{Z}^2$ and
when working with function or distribution spaces we shall freely identify functions/distributions defined on~$\mathbb{T}^2$ \big(resp.\ $\mathbb{R}\times \mathbb{T}^2$\big) with $\mathbb{Z}^2$-periodic functions/distributions defined on $\mathbb{R}^2$ \big(resp.\ $\mathbb{R}\times \mathbb{T}^2$\big).
For $\mathbb{A}$ an~open subset of $\mathbb{T}^2$, $\mathbb{R}\times \mathbb{T}^2$ or $\mathbb{R}^k$ for some $k\geq 1$, we denote by $C(\mathbb{A})$ the space of continuous real valued functions on $\mathbb{A}$ equipped with the supremum norm. As usual for $n\in \mathbb{N}$, $C^{n}$ are the $n$-times continuously differentiable functions. For $\alpha\in \mathbb{R}\setminus \mathbb{N}$, denote by $C^{\alpha}$ the usual H\"older--Besov spaces as in~\cite{Hai14}.

For $n\in \mathbb{N}$, we also use the space of test functions
\[
\mathcal{B}_n:= \{\phi\in C^n(B_1) \mid \|\nabla^m \phi\|_{L^\infty}<\infty,\, m \leq n \},
\]
where $B_1\subset \mathbb{A}$ denotes the unit ball and the choice of $\mathbb{A}$ will be clear from context.

For a function $\rho\colon \mathbb{A}\to \mathbb{R}$, we shall use the following notation for recentering: For $x,y\in \mathbb{A}$, we write $\rho_y(x):= \rho(x-y)$ whenever this is defined.
For any $\lambda \in (0,1]$, we shall use the following scalings:
For $\varphi\colon \mathbb{R}^2\rightarrow \mathbb{R}$, write \smash{$\varphi^\lambda(x)=\lambda^{-2}\varphi\bigl(\frac{x}{\lambda}\bigr)$}, and, similarly, for
$\varphi\colon \mathbb{R}\times \mathbb{R}^2\rightarrow\mathbb{R}$, write \smash{$\varphi^\lambda_x:=\lambda^{-4}\varphi\bigl(\frac{\cdot}{\lambda^2},\frac{\cdot}{\lambda}\bigr)$}.

If there exists a constant $C$ depending only on a parameter set $S$ such that $A\leq CB$, we shall often write this as $A\lesssim_{S} B$.
Finally, we shall use the notation $\fraks$ and $|\fraks|$ from \cite{Hai14}, the latter in our context will be either~$2$ or~$4$.

\section{Solution theory}

In this section, we recall the solution theory for both considered equations. We shall solve the g-PAM equation using regularity structures following~\cite{Hai14} with the necessary modifications for variable coefficients taken from~\cite{Sin25}. For the $\phi^{K+1}_2$-equation, while we could have formulated the fixed point problem similarly using regularity structures, we present it in more elementary language using classical function spaces along the lines of \cite{DD03}.

\subsection{Regularity structures and models}\label{sec:DefRegStructure}
For the purpose here, we define a regularity structure
$(T,G)$ to be a pair, where $T= \bigoplus_{\alpha\in A} T_{\alpha}$ is a finite dimensional normed $A$-graded vector space and $G$ a group acting on $T$ by lower triangular linear maps, i.e., $\gamma \tau - \tau\in T_{<\alpha}:=\bigoplus_{\beta<\alpha}T_\beta$ for every $\tau \in T_\alpha$ and $\gamma\in G$. Throughout, we are only working with one fixed regularity structure.
A model $M=(\Pi,\Gamma)$ consists of a pair of maps: the realisation map $\Pi\colon \mathbb{R}^2\rightarrow {L}\bigl(T,\mathcal{D}'\bigl(\mathbb{R}^d\bigr)\bigr)$
and a map $\Gamma\colon \mathbb{R}^{2}\times \mathbb{R}^2 \to G$ satisfying for every $x,y,z \in \mathbb{R}^{2}$
\[
\Pi_x \circ \Gamma_{xy}= \Pi_y, \qquad \Gamma_{xy} \circ \Gamma_{yz}= \Gamma_{xz},
\]
such that
	\begin{align}\label{eq:mod_bounds_1}
		&\|\Pi\| := \sup_{\alpha \in A}\sup_{\tau \in T_\alpha}\sup_{\varphi \in \mcB_r} \sup_{x \in \mathbb{R}^2} \sup_{\lambda \in (0,1]}
\frac{|\Pi_x \tau (\varphi_x^\lambda) |}{		
		\lambda^{\alpha} |\tau|_\alpha } < \infty,
		\\
		\label{eq:mod_bounds_2}
		&\|\Gamma\| := \sup_{\alpha\in A}\sup_{\tau \in T_\alpha} \sup_{\beta < \alpha } \sup_{x,y \in \mathbb{R}^2}
\frac{|\Gamma_{xy} \tau |_\beta}{	
		 |x-y|_\mfs^{\alpha-\beta} |\tau|_\alpha} < \infty,
	\end{align}
where \smash{$|\cdot |_\beta= | Q_{\beta} (\cdot)|$} with \smash{$Q_\beta\colon T\to T_\beta $} being the canonical projection. We denote by $\mathcal{M}$ the space of all models for the regularity structure $(T,G)$, which we equip with the metric where the distance between two models \smash{$\bigl\{\bigl(\Pi^i, \Gamma^i\bigr)\bigr\}_{1\leq i\leq 2}$} is given by
 replacing $\Pi$, $\Gamma$ in \eqref{eq:mod_bounds_1} and \eqref{eq:mod_bounds_2} by~${\Pi^1 - \Pi^2}$ and $\Gamma^1 - \Gamma^2$, respectively.

\subsection{Solutions to the g-PAM equation}

To define the regularity structure to solve the g-PAM equation, we first fix a set of letters $\{\Xi,\, \partial_i,\, I,\, 1,\, X_i,\, \cdot \mid i=1,2 \}$.
Then, define the set of words
\begin{gather*}
\mathbf{T}= \bigl\{ \Xi,\, \partial_i I\Xi,\, \partial_j I\Xi\cdot \partial_i I\Xi,\, X_i \cdot \Xi,\, \Xi\cdot I\Xi,\, 1,\, I\Xi,\, X_i \mid i,j=1,2\bigr\}
\end{gather*}
and finally set $T$ to be the span of $\mathbf{T}$. In order to define the grading on $T$, we first define homogeneities $| \tau |\in \mathbb{R}$ for each element of $\tau \in \mathbf{T}$ as follows. First set the homogeneities of the letters to be
\[|\Xi|=-1-\kappa,\qquad | \partial_i|=-1, \qquad |I|=2, \qquad |1|=0,\qquad | X_i|=1,\qquad |\cdot|=0\]
 for some $0<\kappa\ll 1$ small enough. The homogeneity of a word is simply the sum of the homogeneities of its letters (counted with multiplicity). Finally, the grading of $T$ is the unique one such that $\tau \in T_{|\tau|}$ for each $\tau\in \mathbf{T}$.
We define the group $G=\bigl(\mathbb{R}^3, + \bigr)$ where the action of $\gamma=(\gamma_{I\Xi}, \gamma_1, \gamma_2 )$ on $T$ is determined by
\begin{itemize}\itemsep=0pt
\item $\gamma \tau = \tau$ for all $\tau \in \bigl\{ \Xi,\, \partial_i I\Xi,\, \partial_j I\Xi\cdot \partial_i I\Xi,\, 1 \mid i,j=1,2\bigr\}$,
\item $\gamma I\Xi= I\Xi + \gamma_{I\Xi}$ and $\gamma X_i= X_i + \gamma_i$ for any $i\in \{1,2\}$,
\item $ \gamma$ is multiplicative with respect to the partial product $\cdot$ on $\mathbf{T}$, i.e.,
$ \gamma (X_i \cdot \Xi) = X_i \cdot \Xi + \gamma_i \Xi$ and
$ \gamma (\Xi\cdot I\Xi)= \Xi \cdot I\Xi + \gamma_{I\Xi} \Xi$.
\end{itemize}

\subsubsection{Canonical and renormalised models}
Recall the definition of the kernel spaces \smash{$\bfK^{2}_{L,R}$} in \cite[Defenition~2.12]{Sin25} (see Appendix~\ref{App:Kernels}) which quantify the regularity assumption \cite[Assumption~5.1]{Hai14}.
For this section, we fix a kernel \smash{$K\in \bfK^{2}_{L,R}$} for some $L,R>1$, which we furthermore assume to be non-anticipative, i.e., for any $x,y\in \mathbb{R}^2$ it holds
$K(t,x,s,y)=0$ whenever $t<s$,
and time-translation invariant,
i.e., $K(t,x,s,y)=K(t-s, x,y)$ for all $t,s\in \mathbb{R}$.
We fix $\kappa\colon \mathbb{R}\to [0,1]$ to be smooth on $(0,\infty)$ such that $\kappa|_{(0,1]}=1$ and $\kappa|_{\mathbb{R}\setminus (0,2]}=0$
 and set for any $x,y\in \mathbb{R}^2$, such that $x\neq y$
\begin{gather}\label{eq:H_kernel}
H(x,y)= \int_{\mathbb{R}} \kappa (s) K(s,x,y)\, \dd s.
\end{gather}
Given a smooth compactly supported function $\rho$ as in \eqref{eq:regularised noise}, we define the associated canonical model $M(\rho)$ to be the pair $(\Pi, \Gamma)$ determined as follows.
The realisation map $\Pi_x\colon T\to \mathcal{D}'\bigl(\mathbb{R}^2\bigr)$ is given by
\begin{gather*}
\Pi_x\Xi(y)= \xi(\rho_y), \qquad \Pi_x 1= 1, \qquad \Pi_x X_i(y)= y_i-x_i, \qquad \Pi_x X_i \cdot \Xi (y)= (y_i-x_i)\xi(\rho_y),
\\
\Pi_x I\Xi(y) = \int_{\mathbb{R}^2} H(y,z)\xi(\rho_z)\,\dd z- \int_{\mathbb{R}^2} H(x,z)\xi(\rho_z)\, \dd z,\\
\Pi_x \partial_i I\Xi(y) =\int_{\mathbb{R}^2} \partial_i H(y,z)\xi(\rho_z)\, \dd z
\end{gather*}
and by
\begin{gather}\label{eq:PAM prod guys}
\Pi_x (\partial_j I\Xi \cdot \partial_i I\Xi) = \Pi_x \partial_j I\Xi \cdot \Pi_x \partial_i I\Xi,
\qquad
\Pi_x (\Xi \cdot I \Xi)= \Pi_x \Xi \cdot \Pi_x I \Xi.
\end{gather}
The map $\Gamma\colon \mathbb{R}^{2} \to G\sim \mathbb{R}^3$ is given by
\[
\Gamma_{xy}= \biggl( \int_{\mathbb{R}^2}(H(x,z)-H(y,z))\xi(\rho_z)\, \dd z,\, x_1-y_1,\,x_2-y_2\biggr).
\]

Given such a canonical model $M=(\Pi, \Gamma)$ and
\smash{$\mathfrak{c}=\bigl( \mathfrak{c}^{\<Xi2>}, \bigl\{\mathfrak{c}^{\<b2>_{ij}}\bigr\}_{i,j=1,2} \bigr)\in C\bigl(\mathbb{T}^2\bigr)^{\times 5}$},
which we~identify
with
continuous periodic functions
 \smash{$\mathfrak{c}^{\<Xi2>}, \mathfrak{c}^{\<b2>_{ij}} \colon \mathbb{R}^2\to \mathbb{R}$}, we define the renormalised model
$\mathfrak{c} \cdot M:= ({\Pi}^\mathfrak{c}, {\Gamma}^\mathfrak{c})$ to be the model that agrees with the model $(\Pi, \Gamma)$ except that for any~${x\in \mathbb{R}^2}$
\begin{gather*}
{\Pi}^\mathfrak{c}_x \partial_j I\Xi \cdot \partial_i I\Xi = {\Pi}_x \partial_j I\Xi \cdot \partial_i I\Xi + \mathfrak{c}^{\<b2>_{ij}},
\qquad
{\Pi}^\mathfrak{c}_x \Xi \cdot I \Xi= {\Pi}_x \Xi \cdot I \Xi +\mathfrak{c}^{\<Xi2>}.
\end{gather*}

Finally, we define the space of admissible models to be the closure of the models described above
\[
\mathcal{M}_{ad}= \overline{\bigl\{ g \cdot M (\rho) \in \mathcal{M} \mid \rho\in C^{\infty}_c(B_{1/2}), g\in C\bigl(\mathbb{T}^2\bigr)^{\times 5} \bigr\}} \subset
 \mathcal{M},
 \]
 see Section \ref{sec:DefRegStructure}.

\subsubsection{Proof of Theorem~\ref{thm:g-pam}}
As a consequence of working with a pathwise solution theory, the only place where the proof of this theorem differs notably from the proof of \cite[Theorem~4.1]{Sin25} is in establishing the necessary stochastic estimates on the model which is the content of Section~\ref{sec:stoch est pam}. In order to make the article relatively self-contained, we present the proof which is relatively compact as it mostly amounts to collecting the necessary ingredients from the existing literature.

\begin{proof}[Proof of Theorem~\ref{thm:g-pam}]
In order to lift our SPDE to an abstract equation on the regularity structure, we fix $\phi\colon \mathbb{R}^2\to [0,1]$ smooth and compactly supported on ${B}_{1+1/10}(0)\subset \mathbb{R}^2$ such that $\sum_{k\in \mathbb{Z}^d} \phi(x+k)=1$ for all $x\in \mathbb{R}^2$ and set for any $t>s$ and $x,y\in \mathbb{R}^2$
\begin{gather*}
K(t,x,s,y)=\kappa(t-s)\phi(x-y)\sum_{k\in \mathbb{Z}^d}\Gamma(t,x,s,y+k),
\end{gather*}
where we recall that $\kappa$ denotes a smooth cut-off from $[0,1]$ to $[0,2]$. Then, it was the content of {\cite[Proposition~2.18]{Sin25}} that for smooth coefficient fields this kernel belongs to the space \smash{$\bfK_{L,R}^2$} for any $L,R>0$, the definition of which we recall in Appendix~\ref{App:Kernels}. But observe that exactly the same argument shows that for $C^{2}$-regular coefficients, $K$ belongs to \smash{$\bfK_{L,R}^2$}, for any $L,R<2$.

Thus, given an admissible model $M\in \mathcal{M}_{ad}$ this allows one to formulate the abstract fixed point problem as in \cite[Section~9]{Hai14}
\begin{gather}\label{abstract pam}
U= \mcK_\gamma \bigl( \mathbf{R}_+ \hat{F}(U, \nabla U)\bigr) + U_{\texttt{in}},
 \end{gather}
 where \[\hat{F}(U, \nabla U):= \sum_{i,j=1}^2 \hat{f}_{i,j}(U)\partial_i U \partial_j U + \hat{g}(U) \Xi. \]

Then, for $L>\gamma> -|\Xi|$ and $\eta\in(0,1)$ \cite[Theorem~2.14]{Sin25}, which is a slight variant of \cite[Theorem~7.8]{Hai14}, shows that there exists $T>0$ and a neighbourhood $\mathcal{O}_M\subset \mathcal{M}_{ad}$ of the model $M$ such that the abstract solution map to \eqref{abstract pam}
\[
S_T\colon\ \mathcal{O}_M \to \mathcal{D}^{\gamma, \eta}_P,\qquad
 M\mapsto U
\] is well defined and continuous.

Next, we consider the sequence of admissible models
\begin{gather}\label{eq:pam sequence of models}
\hat{M}^\delta=\bigl(\hat{\Pi}^\delta,\hat{\Gamma}^{\delta}\bigr):= \mathfrak{c}_\delta \cdot M\bigl(\rho^\delta\bigr), \qquad
\mathfrak{c}_\delta:= \bigl( {c}_{\delta}^{\<Xi2>}, \bigl\{ {c}_{\delta}^{\<b2>_{ij}}\bigr\}_{i,j=1,2} \bigr)
\end{gather}
 with the functions \smash{${c}_{\delta}^{\<Xi2>}$}, \smash{${c}_{\delta}^{\<b2>_{ij}}$} defined in \eqref{eq:renormalisation function pam}. Thus, this proof is complete if one checks the following two points.
\begin{enumerate}\itemsep=0pt
\item For $\delta\in(0,1]$, the function $u_\delta:= \hat{\Pi}^\delta_x \bigl[U^{\delta}(x)\bigr](x)$ for $U^{\delta}= S_T\bigl(\hat{M}^\delta\bigr)$ satisfies \eqref{eq:PAM}.
\item There exists a modification of underlying white noise $\xi$, such that the sequence of models \smash{$\bigl\{\hat{M}^\delta\bigr\}_{\delta>0}$} converge (in probability) to a limiting model \smash{$\hat{M}^0$} as $\delta\to 0$.
\end{enumerate}
The former point is obtained by adapting ad verbatim \cite[Section~9.3]{Hai14}, see also \cite[Section~16.1]{HS23m} or \cite{BB21}. The second point is the content of Proposition~\ref{prop:convergence of the model}.
\end{proof}

\subsection[Solutions to the phi\^{}{K+1}\_2-equation]{Solutions to the $\boldsymbol{\phi^{K+1}_2}$-equation}
In this section we recall the solution theory for the \smash{$\phi^{K+1}_2$}-equation, based on the Da Prato Debussche trick \cite{DD03}.
The mild formulation of the equation \eqref{phi4} with initial condition $u_0$ reads, for any $(t,x)\in (0,\infty)\times \mathbb{R}^2$
\begin{align*}
u_\delta(t,x)={}& \int_{\mathbb{R}^2} \Gamma(t,x,0,y) u_0(y)\,\dd y\\
&{}{+}\, \int_{0}^t \int_{\mathbb{R}^2} \Gamma(t,x,s,y) \bigl(-H_K\bigl(u_\delta(s,y),c_{\delta}^{\<2>}(s,y)\bigr) +\xi_{\delta} (y,s)\bigr)\, \dd s \dd y.
\end{align*}
We define
$\<1>_\delta(t,x):= \int_{\mathbb{R}^{1+2}} \kappa(t-s) \Gamma(t,x,s,y) \xi_{\delta} (s,y)\, \dd s \dd y $.
If we set $v_\delta:= u_\delta-\<1>_\delta$, we find that we can rewrite the integral equation as
\begin{align*}
v_\delta(t,x)+\<1>_\delta(t,x)={}&
 \int_{0}^t \int_{\mathbb{R}^2} \Gamma(t,x,s,y) \bigl(-H_K\bigl(v_\delta(s,y)+\<1>_\delta(s,y),c_{\delta}^{\<2>}(s,y)\bigr) +\xi_{\delta} (y,s)\bigr)\, \dd s\dd y \\
 &{}{+}\,\int_{\mathbb{R}^2} \Gamma(t,x,0,y) u_0(y)\,\dd y.
\end{align*}
Thus using that
$H_{K}(X+Y, C )=\sum_{k=0}^K {K \choose k} X^k H_{K-k}(Y, C)$ and defining
\[
z_\delta (t,y):= - \int_{-\infty}^0 \kappa(t-s) \int_{\mathbb{R}^2} \Gamma\bigl(0,y,s,y'\bigr) \xi_{\delta} \bigl(s,y'\bigr)\, \dd s \dd y',
\]
this can be rewritten as
\begin{align}
v_\delta(t,x)={}& {-}\sum_{k=0}^K {K \choose k} \int_{0}^t \int_{\mathbb{R}^2} \Gamma(t,x,s,y) v^k_\delta(s,y) H_{K-k}\bigl(\<1>_\delta(s,y), c^{\<2>}_\delta(s,y)\bigr)\, \dd s \dd y
\nonumber\\
&{}{+}\,
\int_{\mathbb{R}^2} \Gamma(t,x,0,y) \bigl( u_0(y) + z_\delta(t,y)\bigr)\,\dd y\label{eq:equation after DD}
.
\end{align}
Then, consider for $T>0$ and $\gamma,\beta\in (0,1)$ the weighted space $C_\gamma \bigl((0,T], C^{\beta}\bigl(\mathbb{T}^2\bigr) \bigr)$ of continuous functions $w\colon (0,T]\to C^{\beta}\bigl(\mathbb{T}^2\bigr)$ satisfying
\[
\|w\|_{\gamma,\beta}:=\sup_{0<s\leq T} s^{\gamma} \|w(s)\|_{C^{\beta}} <\infty.
\]
The next proposition is standard, cf.\ \cite[Proposition~4.4]{DD03}, \cite[Theorem~6.2]{MWplane}, \cite[Theorem~3.3]{TW18} or \cite[Proposition~3.2]{HSper}.
\begin{prop}\label{prop:local_wellposedness}
Let $-1/10< \alpha<0$, $\beta,\gamma\in (0,1)$. Let
\[
F=(F_0,\dots,F_K) \in \bigl(C^{\alpha}\bigl((-1,2)\times\mathbb{T}^2\bigr)\bigr)^{K} \qquad \text{and}\qquad w\in C_{\gamma} \bigl((0,2), C^{\beta} \bigr).
\]
Then, there exists $T=T(\|w\|_{\gamma,\beta},\|F\|_{\gamma,\beta}) \in (0,1]$ and a neighborhood
\[
\mathcal{O}_{(w,F)}\subset C_{\gamma} \bigl((0,2], C^{\beta}\bigl(\mathbb{T}^2\bigr) \bigr)\times \bigl(C^{\alpha}\bigl((-1,2)\times\mathbb{T}^2\bigr)\bigr)^{\times K}
\]
of $(w,F)$ such that the solution map
\[S_{T^*}\colon \ \mathcal{O}_{(w,F)} \to C_\gamma((0,T^*], C^{\beta}\bigl(\mathbb{T}^2\bigr), \qquad (w,F)\mapsto v\]
to the equation\footnote{Here we use the standard abuse of notation by writing integrals for the pairing between distributions and functions.}
\[
v(t,x)= w(x,t)
+\sum_{k=0}^K {K \choose k} \int_{0}^t \int_{\mathbb{R}^2} \Gamma(t,x,s,y) v^k(s,y) F_k(s,y)\, \dd s\dd y
\]
is well defined and continuous.
\end{prop}

\begin{proof}[Proof of Theorem~\ref{thm:phi4_2}]
It remains only to check that we can apply
 Proposition~\ref{prop:local_wellposedness} to \eqref{eq:equation after DD},
for example, with $\beta=1/10$, $\gamma=1/5$. Indeed, the map
\[[0,1]\to C_{\gamma} \bigl((0,2], C^{\beta}\bigl(\mathbb{T}^2\bigr)\bigr), \qquad
\delta \mapsto \biggl(t\mapsto\int_{\mathbb{R}^2} \Gamma(t,\cdot,0,y) ( u_0(y) + z_\delta(t,y))\,\dd y\biggr)
\]
is a.s.\ continuous as can be seen by noting that $\mathbb{R}\to C^{\alpha}\bigl(\mathbb{T}^2\bigr)$,
$t\mapsto z_\delta(t,\cdot) $ is continuous for any $\alpha<0$ and then by applying \cite[Lemma~3.1]{HSper}.
Finally, by Lemma~\ref{lem:phi4 estimate} combined with Kolmogorov's continuity criterion
for each $m=0,\dots,K$ and any $\alpha<0$ the map
\[
(0,1]\to C^{\alpha}\bigl((-1,2)\times \mathbb{T}^2\bigr), \qquad \delta\mapsto H_{m}\bigl(\<1>_\delta, c^{\<2>}_\delta\bigr)
\]
extends (up to modification) continuously to $\delta\in [0,1]$.
\end{proof}

\section{Stochastic estimates for the g-PAM equation}\label{sec:stoch est pam}
The aim of this section is to prove the following proposition.
\begin{prop}\label{prop:convergence of the model}
The sequence of models \smash{$\bigl\{\hat{M}^{\delta}\bigr\}_{\delta\in (0,1]}$} defined in \eqref{eq:pam sequence of models} are uniformly bounded and converge, up to a~modification\footnote{Since we use Kolmogorov's continuity criterion several times.} of the underlying noise, to a limiting model $\hat{M}^0$ as ${\delta\to 0}$~a.s.
\end{prop}
This proposition is a well understood consequence of Lemmas~\ref{lem:dumbbell} and~\ref{lem:cherry} below. We present the line of argumentation for the readers convenience.
\begin{proof}

We shall only check the required uniform bounds, the convergence bounds follows analogously. First note that for the polynomial sector there is nothing to check. Thus it remains to check the estimates~\eqref{eq:mod_bounds_1} and~\eqref{eq:mod_bounds_2} (before taking the supremum over $\alpha\in A$ and $\tau\in T_\alpha$) on
\[
 \Xi, \qquad \partial_i I\Xi, \qquad \partial_j I\Xi\cdot \partial_i I\Xi,\qquad X_i \cdot \Xi,\qquad \Xi\cdot I\Xi,\qquad I\Xi.
\]
The almost sure bound on $\hat{ \Pi}^{\delta}_x \Xi= \xi\bigl(\rho^\delta\bigr)$ follows directly from the fact that $\xi \in C^{\alpha}$ a.s. for any $\alpha<-1$, by a standard Kolmogorov criterion. Since, \smash{$\hat{\Gamma}^{\delta}_{xy} \Xi= \Xi$} for all \smash{$x,y\in \mathbb{R}^{2}$}, for the two symbols
$\partial_i I\Xi$ and $ I\Xi$ the estimates~\eqref{eq:mod_bounds_1} and~\eqref{eq:mod_bounds_2} follow from the extension theorem, \cite[Theorem~5.14]{Hai14}. Next, for $\partial_j I\Xi\cdot \partial_i I\Xi$ the bound \eqref{eq:mod_bounds_1} follows again by standard Kolmogorov from Lemma~\ref{lem:cherry} (up to taking a modification of the white noise process) while the bound \eqref{eq:mod_bounds_2} is empty.

It only remains to check the symbols $X_i \cdot \Xi$ and $\Xi\cdot I\Xi$. The bound~\eqref{eq:mod_bounds_2} follows by multiplicativity from the same bound for $X_i $ resp.\ $I\Xi$. Finally, to obtain \eqref{eq:mod_bounds_1} one argues using the countable characterisation of models \cite[Proposition~3.32]{Hai14} (here one takes again a modification of the white noise process).
\end{proof}

In what follows, it will be important to explicitly exhibit the randomness appearing in the kernels in \eqref{eq:frozen kernel} and \eqref{eq:space dep frozen}. We define for $\eta\in \mathbb{R}$
\begin{gather}\label{eq:frozen kernel_factored}
\tilde{Z}^{\eta}(t,y)= \frac{1}{4\pi t \sqrt{\det(A(\eta))}} \exp\biggl(- \frac{ y \cdot A^{-1}(\eta) y}{4t} \biggr),
\\
\label{eq:space dep frozen_factored}
\tilde{G}^{\eta}(y)= \int_{\mathbb{R}} \kappa(t) \tilde{Z}^{\eta}(t,y)\, \dd t, \qquad
\tilde{G}_i^{\eta}(y)= \sum_{j=1}^{2} A^{ij}(\eta) y_j\int_{\mathbb{R}} \kappa(t) \frac{1}{2 t}\tilde{Z}^{\eta}(t,y)\, \dd t,
\end{gather}
since these kernels factor the randomness, i.e., for $h$ as in Assumption~\ref{AssumCoef} it holds that
\[{Z}^{x}(t,y)= \tilde{Z}^{h(x)}(t,y),\qquad {G}^{x}(y)=\tilde{G}^{h(x)}(y), \qquad {G}_i^{x}(y)= \tilde{G}_i^{h(x)}(y).
\]

\subsection[Stochastic estimates on hat{Pi}\_star\^{}{delta} (Xi cdot I*Xi)]{Stochastic estimates on $\boldsymbol{\hat{\Pi}_\star ^\delta (\Xi \cdot I\Xi)}$}
In this subsection, we prove the following lemma.
\begin{Lemma}\label{lem:dumbbell}
For any $\alpha<0$ and $q\in \mathbb{N}$, and $\kappa>0$ small enough,
 \begin{gather*}
 \E\bigl[ \bigl|\hat{\Pi}_\star ^\delta \Xi I\Xi \bigl(\phi^\lambda_\star\bigr)\bigr|^q\bigr]\lesssim_{\alpha,q} \lambda^{\alpha q} \qquad\text{and}\\
 \E\bigl[ \bigl|\hat{\Pi}_\star ^{\delta} \Xi I\Xi \bigl(\phi^\lambda_\star\bigr)-\hat{\Pi}_\star ^{{\delta'}} \Xi I\Xi \bigl(\phi^\lambda_\star\bigr)\bigr|^q\bigr]\lesssim_{\alpha,q} |\delta-{\delta'}|^{\kappa q} \lambda^{(\alpha -\kappa)q},
 \end{gather*}
uniformly over $\star\in \mathbb{R}^2$, $\lambda,\delta,{\delta'}\in (0,1]$ and $\phi \in \mathcal{B}_1$.
\end{Lemma}
\begin{proof}
Unravelling the definition, we find that with $H$ defined in \eqref{eq:H_kernel}
\begin{align*}
\hat{\Pi}_\star^\delta \Xi I\Xi \bigl(\phi^\lambda_\star\bigr)&{}= {\Pi}_\star^\delta \Xi I\Xi \bigl(\phi^\lambda_\star\bigr)
- c^{\<Xi2>}_\delta \bigl(\phi^\lambda_\star\bigr)\\
 &{}= \int_{\mathbb{R}^2} \int_{\mathbb{R}^2}( H(x,y)-H(\star,y))\xi_\delta(x)\xi_\delta(y)\phi^\lambda_\star(x)\, \dd y \dd x
- c^{\<Xi2>}_\delta \bigl(\phi^\lambda_\star\bigr).
\end{align*}

Fix, for $q\in \mathbb{N}$, a smooth function $\psi\colon \mathbb{R}^2\to [0,1]$ compactly supported on ${B}_{1/2q}(0)\subset \mathbb{R}^2$ such that
 \smash{$\psi|_{{B}_{1/4q}(0)}=1$}. Then set for $G^x$ as in \eqref{eq:space dep frozen}
\begin{gather}\label{kernels}
K^{x}(y,z)=\psi(y-z) G^{x}(y-z), \qquad
R(x,y):= H(x,y)-K^{x}(x,y),
\end{gather}
as well as \smash{$\hat{c}_{\delta}^{\<Xi2>}(x) = \int_{\mathbb{R}^2} \psi(y)G^{x}(y) \rho^\delta *\rho^\delta (y)\, \dd y$}.
Thus, we split
\begin{align*}
\hat{\Pi}_\star^\delta \Xi I\Xi (\phi^\lambda_\star)={}& \int_{\mathbb{R}^2}\int_{\mathbb{R}^2} ( H(x,y)-H(\star,y))\xi_\delta(x)\xi_\delta(y) \phi^\lambda_\star(x)\,\dd y \dd x
- \hat{c}^{\<Xi2>}_\delta \bigl(\phi^\lambda_\star\bigr) - \bigl(c^{\<Xi2>}_\delta -\hat{c}^{\<Xi2>}_\delta\bigr) \bigl(\phi^\lambda_\star\bigr)\\
={}& \int_{\mathbb{R}^2}\int_{\mathbb{R}^2} ( K^x(x,y)-K^\star (\star,y)) \xi_\delta(x)\xi_\delta(y) \phi^\lambda_\star(x)\, \dd y \dd x
 \\
 &{}{+} \int_{\mathbb{R}^2}\int_{\mathbb{R}^2} ( R(x,y)-R(\star,y))\xi_\delta(x)\xi_\delta(y) \phi^\lambda_\star(x)\, \dd y \dd x \\
&{}{-}\int_{\mathbb{R}^2}\int_{\mathbb{R}^2} K^x(x,y) \E[ \xi_\delta(x)\xi_\delta(y) ] \phi^\lambda_\star(x)\, \dd y \dd x
- \bigl(c^{\<Xi2>}_\delta -\hat{c}^{\<Xi2>}_\delta\bigr) \bigl(\phi^\lambda_\star\bigr)\\
={}& \int_{\mathbb{R}^2}\int_{\mathbb{R}^2} (K^x(x,y)-K^\star(x,y)) [\xi_\delta(x)\xi_\delta(y)- \E[\xi_\delta(x)\xi_\delta(y)]] \phi^\lambda_\star(x)\, \dd y \dd x\\
&{}{+}\int_{\mathbb{R}^2}\int_{\mathbb{R}^2} (K^\star(x,y)-K^\star(\star,y)) [\xi_\delta(x)\xi_\delta(y)- \E[\xi_\delta(x)\xi_\delta(y)]] \phi^\lambda_\star(x)\, \dd y \dd x\\
&{}{-} \int_{\mathbb{R}^2}\int_{\mathbb{R}^2} K^\star(\star,y) \E[\xi_\delta(x)\xi_\delta(y)] \phi^\lambda_\star(x)\, \dd y \dd x
- \bigl(c^{\<Xi2>}_\delta -\hat{c}^{\<Xi2>}_\delta\bigr) \bigl(\phi^\lambda_\star\bigr)\\
&{}{+} \int_{\mathbb{R}^2}\int_{\mathbb{R}^2} ( R(x,y)-R(\star,y)) \xi_\delta(x)\xi_\delta(y) \phi^\lambda_\star(x)\, \dd y \dd x.
\end{align*}

Estimating the former two lines is the content of Lemmas~\ref{lem:old_term} and \ref{lem:new_term}. The terms on the third line are easily bounded directly. Finally, for the last term, first note that by the classical Levi's construction of the fundamental solution, cf.\ \cite[Proposition~2.6]{Sin25} it follows that
\smash{$\partial_i R\in \bfK^{2}_{L-1,R}$} almost surely for any $L\in (1,2)$ and that the random variable $\|\partial_i R\|_{\beta; L-1,R}$ has a finite $q$-th moment for any $q<\infty$. Since, for any $\eps>0$ the random variable $ \sup_{\delta\in [0,1]} \| \xi_\delta\|_{C^{-1-\eps}}$ has also finite $q$-th moment, it follows that for any $q<\infty$
\begin{equation}\label{EstimateErrorTermKernel}
\sup_{\delta\in (0,1)} \biggl\|\int_{\mathbb{R}^2} R(\cdot,y) \xi_\delta(y)\, \dd y\biggr\|_{C^{L}} \in L^{q}(\Omega).
\end{equation}
Essentially\footnote{This is for example exactly the reconstruction bound, \cite[equation~(3.3)]{Hai14} when proving of Young multiplication using regularity structures \cite[Proposition~4.14]{Hai14}.} by Young multiplication,
it follows that for any $\kappa>0$
\[ \sup_{\lambda\in (0,1]} \lambda^{-\kappa}\sup_{\delta\in [0,1]} \biggl|\int_{\mathbb{R}^2}\int_{\mathbb{R}^2} ( R(x,y)-R(\star,y)) \xi_\delta(x)\xi_\delta(y) \phi^\lambda_\star(x)\, \dd y \dd x\biggr|\in L^{q}(\Omega).\]
The difference bound for $\delta, \delta' \in (0,1]$ follows similarly.
\end{proof}

In the next lemma, we use the notation $\xi_\delta(x)\diamond \xi_\delta(y):= \xi_\delta(x)\xi_\delta(y)- \E[\xi_\delta(x)\xi_\delta(y)],$
cf.~Appendix~\ref{App:Gaussian int}.
\begin{Lemma}\label{lem:old_term}
 For any $\alpha<0$ and $q\in 2\mathbb{N}$, in the setting of the proof of Lemma~{\rm \ref{lem:dumbbell}}, it holds that
\begin{align}
&\E \biggl[ \biggl|\int_{\mathbb{R}^2}\int_{\mathbb{R}^2} (K^\star(x,y)-K^\star(\star,y)) [ \xi_\delta(x)\diamond \xi_\delta(y)] \phi^\lambda_\star(x)\,\dd y \dd x \biggr|^q \biggr]^{\frac{1}{q}}\lesssim_{q,\alpha} \lambda^\alpha, \label{eq:bond...}
\end{align}
as well as
\begin{align}
&\E \biggr[ \biggr|\int_{\mathbb{R}^2}\int_{\mathbb{R}^2} (K^\star(x,y)-K^\star(\star,y)) [\xi_\delta(x)\diamond \xi_\delta(y) - \xi_{\delta'}(x)\diamond \xi_{\delta'}(y) ]] \phi^\lambda_\star(x)\, \dd y \dd x \biggr|^q \biggr]^{\frac{1}{q}} \nonumber\\
&\qquad{}\lesssim_{q,\alpha} \bigl|\delta-{\delta'}\bigr|^\kappa \lambda^{\alpha-\kappa}, \label{eq:diffbond...}
\end{align}
uniformly over $\star\in \mathbb{R}^2$, $\lambda,\delta, \delta' \in (0,1]$ and $\phi \in \mathcal{B}_1$.
 \end{Lemma}
 \begin{proof}
 For the proof, we factor the randomness in the kernel by defining
\[\tilde{K}^{\eta}(y,z)=\psi(y-z) \tilde{G}^{\eta}(y-z), \]
for $\tilde{G}$ as in \eqref{eq:space dep frozen_factored}.
We shall also use the suggestive notation $\partial_{\uparrow}$ to denote differentiation in the `upper' $\eta$-variable.
 Furthermore, assume without loss of generality that \smash{$\lambda<\frac{1}{2q}$}. Then, expanding the $q$-th power and using the Gaussian integration by parts Lemma \ref{lem:wick_integration by parts}, the term to bound is
 \begin{align}
&\E \biggl[ \biggl| \int_{\mathbb{R}^2}\int_{\mathbb{R}^2} \bigl(\tilde{K}^{h(\star)}(x,y)-\tilde{K}^{h(\star)}(\star,y)\bigr) [\xi_\delta(x)\xi_\delta(y)-\mathbb{E}[\xi_\delta(x)\xi_\delta(y)]] \phi^\lambda_\star(x)\,\dd y \dd x \biggr|^q \biggr] \nonumber\\
&\qquad{}=\sum_{
{\underset{\vert J\vert\text{ even}}
{J\subset\{1,\dots,2q\}}}}
\sum_{P\in \mathcal{P}_2(J)} \int \Biggl[\prod_{i=1}^q \phi^\lambda_{\star}(x_{2i-1})\prod_{\{i,j\}\in P}\mathbb{E}[\xi_\delta(x_i)\xi_\delta(x_j)] \nonumber \\
&\qquad\qquad{} \times\mathbb{E}\Biggl[\partial_{\uparrow}^{\vert J^c\vert}
\prod_{i=1}^{q}\bigl(\tilde{K}^{h(\star)}(x_{2i-1},x_{2i})-\tilde{K}^{h(\star)}(\star,x_{2i})\bigr)
\Biggr]\prod_{i\in J^c}\mathbb{E}[\xi_\delta(x_i) \xi * \sigma(\star)] \Biggr], \label{eq:rewrite}
\end{align}
where the integral is over \smash{$(x_1,\dots,x_{2q}) \in \bigl(\mathbb{R}^{2}\bigr)^{2q}$}. Inserting
\begin{gather}\label{eq:correlation identities}
\mathbb{E}[\xi_\delta(x_i)\xi_\delta(x_j)]= \bigl(\rho^{*2}\bigr)^\delta(x_i-x_j),
\qquad
F^{\star}_\delta(x):=\mathbb{E}[\xi_\delta(x) \xi * \sigma(\star)]
\end{gather}
and using the binomial formula for derivatives of products, we can bound \eqref{eq:rewrite} by
\begin{align*}
&
\sum_{{\underset{\vert J\vert\text{ even}}
{J\subset\{1,\dots,2q\}}}}
\sum_{P\in \mathcal{P}_2(J)}
\sum_{
\underset{|k|= |J^c|}
{k\in \mathbb{N}^q }}
\Biggl|\int \prod_{i=1}^q \phi^\lambda_{\star}(x_{2i-1})\prod_{\{i,j\}\in P}\bigl(\rho^{*2}\bigr)^\delta(x_i-x_j)\\
&\hphantom{\sum_{{\underset{\vert J\vert\text{ even}}{J\subset\{1,\dots,2q\}}}}\sum_{P\in \mathcal{P}_2(J)}\sum_{\underset{|k|= |J^c|}{k\in \mathbb{N}^q }}\Biggl|\int}{}
\times \mathbb{E}\Biggl[
\prod_{i=1}^{q}\partial_{\uparrow}^{k_i} \bigl(\tilde{K}^{h(\star)}(x_{2i-1},x_{2i})-\tilde{K}^{h(\star)}(\star,x_{2i})\bigr)
\Biggr]
\prod_{i\in J^c}
F^{\star}_{\delta}(x_i)\Biggr|,
\end{align*}
 each summand of which can in turn be bounded by
\begin{gather*}
\mathbb{E}\Biggl[\Biggl| \int {\prod_{i=1}^q \phi^\lambda_{\star}(x_{2i-1})}\prod_{\{i,j\}\in P}\bigl(\rho^{*2}\bigr)^\delta(x_i-x_j) \\
\hphantom{\mathbb{E}\Biggl[\Biggl|\int }{}
\times
\prod_{i=1}^{q}\partial_{\uparrow}^{k_i} \bigl(\tilde{K}^{h(\star)}(x_{2i-1},x_{2i})-\tilde{K}^{h(\star)}(\star,x_{2i})\bigr)
\prod_{i\in J^c}
F^{\star}_{\delta}(x_i) \Biggr| \Biggr].
\end{gather*}
This integral can be graphically represented using a directed graph $\CCG=(\CCV,\CCE)$ constructed as follows.
The vertex set is given by $\CCV =\{v_0,\dots,v_{2q}\}$, each element $v_i$ for $i\geq 1$ of which we interpret as representing the integration variable $x_i$ and $v_0$ representing $\star$. We shall apply the Hairer--Quastel criterion, Theorem~\ref{thm:HQ}, and therefore use notations from Appendix~\ref{App:HairerQuastel}.
 \begin{itemize}\itemsep=0pt
\item The graph contains $q$ distinguished edges representing the generic test function $\phi_\star^\lambda$, which we shall illustratively draw as
 \tikz[baseline=-0.1cm] \draw[testfcn] (1,0) to (0.0,0) \enlarge; in the schematics below, one directed from $v_0$ to each odd vertex $v_{2i-1}$.
For each such edge $e$, we set $(a_e,r_e)=(0, 0)$.
\item It contains $q$ edges, for each $i=1,\dots,q$ one representing a factors of $\partial_{\uparrow}^{k_i}\bigl(\tilde{K}^{h(\star)}(x_{2i-1},x_{2i})-\tilde{K}^{h(\star)}(\star,x_{2i})\bigr)$ directed from the vertex $v_{2i}$ to $v_{2i-1}$.
We draw each such edge $e$ as \tikz[baseline=-0.1cm] \draw[kernel1] (0.0,0) to (1,0); and set $(a_e, r_e)=(\kappa',1)$ for some $\kappa'\in (0,1)$ which we specify later. We call these $K$-type edges.
\item It contains $|J|/2$ edges representing \smash{$\bigl(\rho^{*2}\bigr)^\delta$}, drawn as \tikz[baseline=-0.1cm] \draw[rho] (0.0,0) to (1,0); (we suppress the arbitrary orientation since it is irrelevant). For each $P\in \mathcal{P}_2(J)$ and $\{i,j\} \in P$ one connecting $v_i$ to $v_j$.
 For such an edge $e$, we set $I_e=1$ and $(a_e,r_e)=(2+\kappa, -1)$ for some $\kappa\in (0,1)$ specified later. We call these $\rho$-type edges.
 \item Finally, there are $|J^c|$ edges, for each $i \in J^c$ one representing the function $F^\star_\delta(x_i)$ which connects $\star$ to $v_i$. For such an edge $e$, we set $(a_e, r_e)=(0,0)$.
 We shall omit these edges in the schematic drawings since they
represent functions that are uniformly bounded in $C^1$ and thus do not contribute to the argument used to bound the singular integral.
 \end{itemize}

Let us illustrate this by listing (up to symmetries and the functions $F^\star$ omitted in the drawing) all possibilities for the case $q=2$
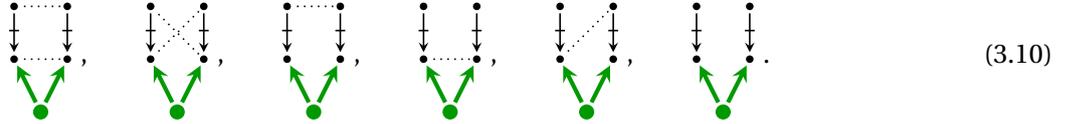
\begin{equation}\label{eq:illustration p=2}
\begin{tikzpicture}[scale=0.35,baseline=0.3cm]
	\node at (0,-1) [root] (root) {};
	\node at (-1,1) [dot] (left) {};
	\node at (-1,3) [dot] (left1) {};
	\node at (1,1) [dot] (variable1) {};
	\node at (1,3) [dot] (variable2) {};
	
	\draw[testfcn] (left) to (root);
	\draw[testfcn] (variable1) to (root);
	
	\draw[kernel1] (left1) to (left);
	\draw[kernel1] (variable2) to (variable1);
	\draw[rho] (variable2) to (left1);
	\draw[rho] (variable1) to (left);
\end{tikzpicture},
\qquad
\begin{tikzpicture}[scale=0.35,baseline=0.3cm]
	\node at (0,-1) [root] (root) {};
	\node at (-1,1) [dot] (left) {};
	\node at (-1,3) [dot] (left1) {};
	\node at (1,1) [dot] (variable1) {};
	\node at (1,3) [dot] (variable2) {};
	
	\draw[testfcn] (left) to (root);
	\draw[testfcn] (variable1) to (root);
	
	\draw[kernel1] (left1) to (left);
	\draw[kernel1] (variable2) to (variable1);
	\draw[rho] (variable2) to (left);
	\draw[rho] (variable1) to (left1);
\end{tikzpicture},
\qquad
\begin{tikzpicture}[scale=0.35,baseline=0.3cm]
	\node at (0,-1) [root] (root) {};
	\node at (-1,1) [dot] (left) {};
	\node at (-1,3) [dot] (left1) {};
	\node at (1,1) [dot] (variable1) {};
	\node at (1,3) [dot] (variable2) {};
	
	\draw[testfcn] (left) to (root);
	\draw[testfcn] (variable1) to (root);
	
	\draw[kernel1] (left1) to (left);
	\draw[kernel1] (variable2) to (variable1);
	\draw[rho] (variable2) to (left1);
\end{tikzpicture},
\qquad
\begin{tikzpicture}[scale=0.35,baseline=0.3cm]
	\node at (0,-1) [root] (root) {};
	\node at (-1,1) [dot] (left) {};
	\node at (-1,3) [dot] (left1) {};
	\node at (1,1) [dot] (variable1) {};
	\node at (1,3) [dot] (variable2) {};
	
	\draw[testfcn] (left) to (root);
	\draw[testfcn] (variable1) to (root);
	
	\draw[kernel1] (left1) to (left);
	\draw[kernel1] (variable2) to (variable1);
	\draw[rho] (variable1) to (left);
\end{tikzpicture},
\qquad
\begin{tikzpicture}[scale=0.35,baseline=0.3cm]
	\node at (0,-1) [root] (root) {};
	\node at (-1,1) [dot] (left) {};
	\node at (-1,3) [dot] (left1) {};
	\node at (1,1) [dot] (variable1) {};
	\node at (1,3) [dot] (variable2) {};
	
	\draw[testfcn] (left) to (root);
	\draw[testfcn] (variable1) to (root);
	
	\draw[kernel1] (left1) to (left);
	\draw[kernel1] (variable2) to (variable1);
	\draw[rho] (variable2) to (left);
\end{tikzpicture},
\qquad
\begin{tikzpicture}[scale=0.35,baseline=0.3cm]
	\node at (0,-1) [root] (root) {};
	\node at (-1,1) [dot] (left) {};
	\node at (-1,3) [dot] (left1) {};
	\node at (1,1) [dot] (variable1) {};
	\node at (1,3) [dot] (variable2) {};
	
	\draw[testfcn] (left) to (root);
	\draw[testfcn] (variable1) to (root);
	
	\draw[kernel1] (left1) to (left);
	\draw[kernel1] (variable2) to (variable1);
\end{tikzpicture}.
\end{equation}
Let us also provide some examples for the case $q=5$
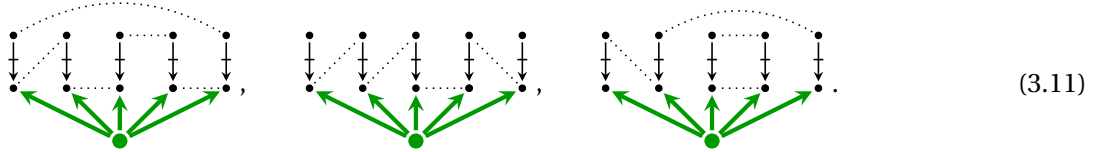
\begin{equation}\label{eq:illustration p=5}
\begin{tikzpicture}[scale=0.35,baseline=0.3cm]
	\node at (0,-1) [root] (root) {};
	\node at (-2,1) [dot] (ld) {};
	\node at (-4,1) [dot] (lld) {};
	\node at (-2,3) [dot] (lu) {};
	\node at (-4,3) [dot] (llu) {};
	\node at (0,3) [dot] (mu) {};
	\node at (2,1) [dot] (rd) {};
	\node at (4,1) [dot] (rrd) {};
	\node at (2,3) [dot] (ru) {};
	\node at (2,1) [dot] (md) {};
	\node at (4,3) [dot] (rru) {};
	\node at (0,1) [dot] (md) {};	
	
	
	\draw[testfcn] (ld) to (root);
	\draw[testfcn] (lld) to (root);
	\draw[testfcn] (md) to (root);
	\draw[testfcn] (rd) to (root);
	\draw[testfcn] (rrd) to (root);
	
	\draw[kernel1] (lu) to (ld);
	\draw[kernel1] (llu) to (lld);
	\draw[kernel1] (mu) to (md);
	\draw[kernel1] (ru) to (rd);
	\draw[kernel1] (rru) to (rrd);
	\draw[rho] (ru) to (mu);
	\draw[rho] (rrd) to (rd);
	\draw[rho] (ld) to (md);
	\draw[rho] (lld) to (lu);
	\draw[rho] (rru) to[bend right=30] (llu);
\end{tikzpicture},
\qquad
\begin{tikzpicture}[scale=0.35,baseline=0.3cm]
	\node at (0,-1) [root] (root) {};
	\node at (-2,1) [dot] (ld) {};
	\node at (-4,1) [dot] (lld) {};
	\node at (-2,3) [dot] (lu) {};
	\node at (-4,3) [dot] (llu) {};
	\node at (0,3) [dot] (mu) {};
	\node at (2,1) [dot] (rd) {};
	\node at (4,1) [dot] (rrd) {};
	\node at (2,3) [dot] (ru) {};
	\node at (2,1) [dot] (md) {};
	\node at (4,3) [dot] (rru) {};
	\node at (0,1) [dot] (md) {};	
	
	
	\draw[testfcn] (ld) to (root);
	\draw[testfcn] (lld) to (root);
	\draw[testfcn] (md) to (root);
	\draw[testfcn] (rd) to (root);
	\draw[testfcn] (rrd) to (root);
	
	\draw[kernel1] (lu) to (ld);
	\draw[kernel1] (llu) to (lld);
	\draw[kernel1] (mu) to (md);
	\draw[kernel1] (ru) to (rd);
	\draw[kernel1] (rru) to (rrd);
	\draw[rho] (rd) to (md);
	\draw[rho] (rrd) to (ru);
	\draw[rho] (ld) to (mu);
	\draw[rho] (lld) to (lu);
\end{tikzpicture},
\qquad
\begin{tikzpicture}[scale=0.35,baseline=0.3cm]
	\node at (0,-1) [root] (root) {};
	\node at (-2,1) [dot] (ld) {};
	\node at (-4,1) [dot] (lld) {};
	\node at (-2,3) [dot] (lu) {};
	\node at (-4,3) [dot] (llu) {};
	\node at (0,3) [dot] (mu) {};
	\node at (2,1) [dot] (rd) {};
	\node at (4,1) [dot] (rrd) {};
	\node at (2,3) [dot] (ru) {};
	\node at (2,1) [dot] (md) {};
	\node at (4,3) [dot] (rru) {};
	\node at (0,1) [dot] (md) {};	
	
	
	\draw[testfcn] (ld) to (root);
	\draw[testfcn] (lld) to (root);
	\draw[testfcn] (md) to (root);
	\draw[testfcn] (rd) to (root);
	\draw[testfcn] (rrd) to (root);
	
	\draw[kernel1] (lu) to (ld);
	\draw[kernel1] (llu) to (lld);
	\draw[kernel1] (mu) to (md);
	\draw[kernel1] (ru) to (rd);
	\draw[kernel1] (rru) to (rrd);
	\draw[rho] (ru) to (mu);
	\draw[rho] (rd) to (md);
	\draw[rho] (llu) to (ld);
	\draw[rho] (rru) to[bend right=30] (lu);
\end{tikzpicture}.
\end{equation}
We shall check the Hairer--Quastel criterion for the graphs obtained after adding `fictitious' edges connecting the remaining vertices $\{v_{i}\}_{i\geq 1}$ which do not already have a $\rho$-type incident edge (in~some arbitrary way).
These new edges represent simply the constant $1$, but we artificially set $(a_e,r_e)=({1+\kappa''},0)$ for some $\kappa''\in (0,1)$ specified later.
Below we illustrate using the (last four) diagrams of
 \eqref{eq:illustration p=2} and the middle diagram in \eqref{eq:illustration p=5} and
 drawing the new edges as 
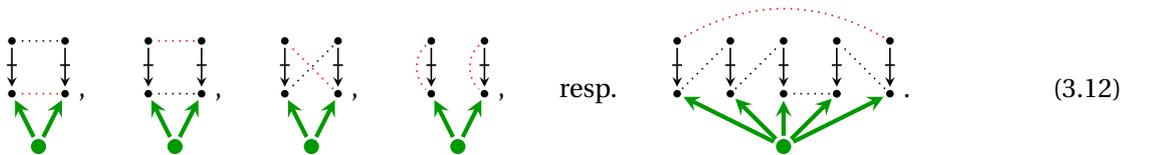
\tikz[baseline=-0.1cm] \draw[red, rho] (0.0,0) to (1,0);
\begin{equation*}
\begin{tikzpicture}[scale=0.35,baseline=0.3cm]
	\node at (0,-1) [root] (root) {};
	\node at (-1,1) [dot] (left) {};
	\node at (-1,3) [dot] (left1) {};
	\node at (1,1) [dot] (variable1) {};
	\node at (1,3) [dot] (variable2) {};
	
	\draw[testfcn] (left) to (root);
	\draw[testfcn] (variable1) to (root);
	
	\draw[kernel1] (left1) to (left);
	\draw[kernel1] (variable2) to (variable1);
	\draw[rho] (variable2) to (left1);
		\draw[red, rho] (variable1) to (left);
\end{tikzpicture},
\qquad
\begin{tikzpicture}[scale=0.35,baseline=0.3cm]
	\node at (0,-1) [root] (root) {};
	\node at (-1,1) [dot] (left) {};
	\node at (-1,3) [dot] (left1) {};
	\node at (1,1) [dot] (variable1) {};
	\node at (1,3) [dot] (variable2) {};
	
	\draw[testfcn] (left) to (root);
	\draw[testfcn] (variable1) to (root);
	
	\draw[kernel1] (left1) to (left);
	\draw[kernel1] (variable2) to (variable1);
	\draw[red, rho] (variable2) to (left1);
	\draw[rho] (variable1) to (left);
\end{tikzpicture},
\qquad
\begin{tikzpicture}[scale=0.35,baseline=0.3cm]
	\node at (0,-1) [root] (root) {};
	\node at (-1,1) [dot] (left) {};
	\node at (-1,3) [dot] (left1) {};
	\node at (1,1) [dot] (variable1) {};
	\node at (1,3) [dot] (variable2) {};
	
	\draw[testfcn] (left) to (root);
	\draw[testfcn] (variable1) to (root);
	
	\draw[kernel1] (left1) to (left);
	\draw[kernel1] (variable2) to (variable1);
	\draw[rho] (variable2) to (left);
	\draw[red, rho] (variable1) to (left1);
\end{tikzpicture},
\qquad
\begin{tikzpicture}[scale=0.35,baseline=0.3cm]
	\node at (0,-1) [root] (root) {};
	\node at (-1,1) [dot] (left) {};
	\node at (-1,3) [dot] (left1) {};
	\node at (1,1) [dot] (variable1) {};
	\node at (1,3) [dot] (variable2) {};
	
	\draw[testfcn] (left) to (root);
	\draw[testfcn] (variable1) to (root);
	
	\draw[kernel1] (left1) to (left);
	\draw[kernel1] (variable2) to (variable1);
	\draw[red, rho] (left1) to [bend right=50] (left);
	\draw[red, rho] (variable2) to [bend right=50] (variable1);
\end{tikzpicture}, \qquad \text{resp.} \quad
\begin{tikzpicture}[scale=0.35,baseline=0.3cm]
	\node at (0,-1) [root] (root) {};
	\node at (-2,1) [dot] (ld) {};
	\node at (-4,1) [dot] (lld) {};
	\node at (-2,3) [dot] (lu) {};
	\node at (-4,3) [dot] (llu) {};
	\node at (0,3) [dot] (mu) {};
	\node at (2,1) [dot] (rd) {};
	\node at (4,1) [dot] (rrd) {};
	\node at (2,3) [dot] (ru) {};
	\node at (2,1) [dot] (md) {};
	\node at (4,3) [dot] (rru) {};
	\node at (0,1) [dot] (md) {};	
	
	
	\draw[testfcn] (ld) to (root);
	\draw[testfcn] (lld) to (root);
	\draw[testfcn] (md) to (root);
	\draw[testfcn] (rd) to (root);
	\draw[testfcn] (rrd) to (root);
	
	\draw[kernel1] (lu) to (ld);
	\draw[kernel1] (llu) to (lld);
	\draw[kernel1] (mu) to (md);
	\draw[kernel1] (ru) to (rd);
	\draw[kernel1] (rru) to (rrd);
	\draw[rho] (rd) to (md);
	\draw[rho] (rrd) to (ru);
	\draw[rho] (ld) to (mu);
	\draw[rho] (lld) to (lu);
	
	\draw[red, rho] (rru) to[bend right=30] (llu);
\end{tikzpicture}.
\end{equation*}

Before checking the assumptions of the Hairer--Quastel criterion, let us make the following observations about any diagram appearing in the above sum.
Consider the (multi) graph obtained after removing the distinguished edges (illustrated above in green) and the edges representing the functions $F_\delta^{\star}$ (which are already omitted in the drawing).
Every vertex of this new graph has exactly two incident edges, one of $K$-type and one which is either of $\rho$-type or fictitious. Thus, its connected components are either cycle graphs or a multi graphs with exactly to vertices and two edges one of which fictitious.
Next, we check the Assumptions~\ref{ass:mainGraph}.
\begin{enumerate}\itemsep=0pt
\item The first condition that $a_{e}+r_e<2$ is clearly satisfied if $\kappa'+ \kappa''<1$ (where
we recall that for a multi-edge we simply add the corresponding labels).
\item
For the second condition, since the left-hand side of \eqref{e:assEdges}
only involves internal edges, we can assume that \smash{$\bar \CCV= \bigcup_{e\in\bar{\CCE}} e$} for some subset $\bar{\CCE}\subset \CCE$.
We first check the inequality for connected components of the subgraph $\bar{\CCG}=(\CCV, \CCE)$ with at least $3$ vertices, where we have two cases:
\begin{enumerate}\itemsep=0pt
\item If $\bar \CCG$ is a cycle graph, it contains $2M$ edges, half of which are of $\rho$-type or fictitious and half of $K$-type for $M\geq 2$. Thus
\[ \sum_{e \in \CCE_0 (\bar \CCV)} a_e= M(2+ \kappa + \kappa') <4M-2 = 2 \bigl(|\bar \CCV\| - 1\bigr).\]
\item
If it is not a cycle, but has an even number of edges $2M$ for $M\geq 1$, we find as above
\[ \sum_{e \in \CCE_0 (\bar \CCV)} a_e =M(2+ \kappa + \kappa') < 4M = 2\bigl(|\bar \CCV| - 1\bigr).\]
If it has an odd number $2M+1$ of edges there is one more of either type of edge and whenever $M\geq 1$
\[\sum_{e \in \CCE_0 (\bar \CCV)} a_e \leq M(2+ \kappa + \kappa') + 2+\kappa
< 4M +2 =2\bigl(|\bar \CCV| - 1\bigr).\]
\end{enumerate}
Thus it remains to check the inequality for sub-graphs of $\CCG$ where some connected components might consist of only two vertices. Since $\vert \bar{\CCV}\vert\geq 3$ by assumption, this can not be the whole graph. Thus, we conclude by noting that
 the inequality holds when the graph consists of only isolated (multi) edges or of isolated (multi) edges and connected components which already satisfy the inequality by themselves.

\item For the third point, note that our setting \eqref{e:assEdges2} simplifies to checking for any $\bar \CCV \subset \CCV_0$ of cardinality at least $1$ that
\[
\sum_{e \in \CCE_0 (\bar \CCV) } a_e
+\sum_{e \in \CCE^{\uparrow} (\bar \CCV) } \kappa'
 - \bigl| \CCE^{\downarrow} (\bar \CCV) \bigr| < 2|\bar \CCV|.
\]
We consider the graph $\bar{\CCG}\subset \CCG$ with vertex set $\bar{\CCV}$ and maximal number of edges.
As in the previous case, we decompose it into connected components and note that by additivity it suffices to check the inequality for each of those.
For singletons, it clearly holds.
For graphs consisting of only one edge, one again simply checks all possible cases.
Finally, for higher cardinality the claim follows from the former item.

\item To check the fourth point of Assumption~\ref{ass:mainGraph}, since every element of $\bar{\CCV}$ has an adjacent outgoing $K$-type edge\footnote{That is, it is drawn `at the top' in the illustrations.} and the set
$\CCE^{\downarrow}(\bar\CCV)$ is empty, \eqref{e:assEdges3} simplifies to
\[
 \sum_{e\in \CCE(\bar\CCV) } a_e
 +\bigl|\CCE^{\uparrow}(\bar\CCV) \bigr|
> 2 |\bar \CCV|.
\]
In turn, since each of these top edges has exactly one outgoing edge with $r_e>0$, this reduces to \smash{$\sum_{e\in \CCE(\bar\CCV) } a_e
> |\CCV |$} which is true since there is either an incident edge of type $\rho$ and if not a fictitious one.
\end{enumerate}
Thus we complete the proof of \eqref{eq:bond...} by noting that for any $\alpha<0$ we can choose $\kappa$, $\kappa'$ such that
\[\alpha q \gtrsim 2 |\CCV\setminus \CCV_\star| - \sum_{e\in \CCE} a_e \geq 2q - q\kappa'- q (2+ \kappa) = -q(\kappa + \kappa') \]
and all the kernels have norms with finite moments.
Finally, to obtain the difference bound \eqref{eq:diffbond...} one argues exactly the same way, but replacing in a telescopic sum each occurrence of
\smash{$\bigl(\rho^{*2}\bigr)^{\delta}$} by \smash{$\bigl(\rho^{*2}\bigr)^{\delta}-\bigl(\rho^{*2}\bigr)^{{\delta'}}$} and
of the function $F^{\star}_\delta$ by
$
F^{\star}_{\delta,{\delta'}}(x):= \mathbb{E}[(\xi_{\delta}- \xi_{{\delta'}})(x) (\xi * \sigma)(\star)].
$
\end{proof}

\begin{Lemma}\label{lem:new_term}
For $\alpha<0$, $q\in 2\mathbb{N}$ and $\kappa\in (0,1)$ in the setting of the proof of Lemma~{\rm \ref{lem:dumbbell}}, it holds that
\begin{gather}\label{to bound new term}
\E \Biggl[ \Biggl|\int_{\mathbb{R}^2}\int_{\mathbb{R}^2} (K^x(x,y)-K^\star(x,y)) [\xi_\delta(x)\diamond \xi_\delta(y)] \phi^\lambda_\star(x)\,\dd y \dd x\Biggr|^q \Biggr]^{\frac{1}{q}}\lesssim_{q,\alpha} \lambda^\alpha,
\end{gather}
as well as
\begin{align}
&\E \Biggl[ \Biggl|\int_{\mathbb{R}^2}\int_{\mathbb{R}^2} (K^x(x,y)-K^\star(x,y)) [\xi_{\delta}(x)\diamond\xi_{\delta}(y)-\xi_{{\delta'}}(x)\diamond\xi_{{\delta'}}(y) ] \phi^\lambda_\star(x)\,\dd y \dd x \Biggr|^q \Biggr]^{\frac{1}{q}} \nonumber\\
&\qquad
\lesssim_{q,\alpha} |\delta-{\delta'}|^\kappa \lambda^{\alpha-\kappa}, \label{eq:diffbond2...}
\end{align}
uniformly over $\star\in \mathbb{R}^2$, $\lambda,\delta,\delta'\in (0,1]$ and $\phi \in \mathcal{B}_1$.
\end{Lemma}
\begin{proof}
We write for any $x,y\in \mathbb{R}^2$, \smash{$\tilde{K}^{h(x)}(x,y)-\tilde{K}^{h(\star)}(x,y) = \sum_{i=1,2} (x_i-\star_i) \tilde{K}^{(i)}(x,y)$}, where
\[\tilde{K}^{(i)}(x,y):=
 \int_{0}^1 \partial_{i} h(s x+ (1-s)\star) \cdot (\partial_{\uparrow} K)^{h(\lambda x+ (1-\lambda)\star) } (x,y)\, \dd s
.\]
Defining the modified test-function
\smash{$\bigl(\tilde{\phi}^{(i)}\bigr)^\lambda\hspace{-0.3pt} (x)\!:=\! \frac{x_i}{\lambda} \phi^\lambda(x)$},
we bound the left-hand side of~\eqref{to bound new term}~by
\[
\sum_{i=1}^2 \lambda\E \Biggl[ \Biggl| \int \tilde{K}^{(i)}(x,y)[\xi_\delta(x)\xi_\delta(y)- \E[\xi_\delta(x)\xi_\delta(y)]] \bigl(\tilde{\phi}^{(i)}\bigr)^\lambda_\star(x)\Biggr|^q \Biggr]^{1/q}.
\]
Thus, we are left to prove that
\[\E \Biggr[ \Biggl| \int \tilde{K}^{(i)}(x,y)[\xi_\delta(x)\xi_\delta(y)- \E[\xi_\delta(x)\xi_\delta(y)]] \bigl(\tilde{\phi}^{(i)}\bigr)^\lambda_\star(x)\Biggr|^q \Biggr] \lesssim \lambda^{(\alpha-1)q}\]
for which one argues as in the proof of Lemma~\ref{lem:old_term}, but employing the slight variant of the Hairer--Quastel criterion, Theorem~\ref{thm:slight variant of HQ} and with the additional observation that \smash{$\|\nabla h\|_{L^{\infty}(\mathbb{T}^2)}\in L^{q}(\Omega)$} by Assumption~\ref{AssumCoef}.
Finally, \eqref{eq:diffbond2...} is obtained by an analogous argument.
 \end{proof}

\subsection[Stochastic estimates on hat{Pi}\_star\^{}{delta} (partial\_i I Xi cdot partial\_j I Xi)]{Stochastic estimates on $\boldsymbol{\hat{\Pi}_\star ^\delta (\partial_i I \Xi \cdot \partial_j I\Xi)}$}
In this subsection, we prove the following lemma.
\begin{Lemma}\label{lem:dumbbell2}
For any $\alpha<0$ and $q\in \mathbb{N}$, and $\kappa>0$ small enough,
\begin{gather*}
\E\bigl[ \bigl|\hat{\Pi}_\star ^\delta \partial_i I \Xi \cdot \partial_j I\Xi \bigl(\phi^\lambda_\star\bigr)\bigr|^q\bigr]\lesssim_{q,\alpha} \lambda^{q \alpha}
\qquad \text{and} \\
\E\bigl[ \bigl|\hat{\Pi}_\star ^\delta \partial_i I \Xi \cdot \partial_j I\Xi \bigl(\phi^\lambda_\star\bigr) -\hat{\Pi}_\star ^{\delta'} \partial_i I \Xi \cdot \partial_j I\Xi \bigl(\phi^\lambda_\star\bigr)\bigr|^q
\bigr]\lesssim_{\alpha,q} |\delta-{\delta'}|^{\kappa q} \lambda^{(\alpha -\kappa)q},
\end{gather*}
uniformly over $\star\in \mathbb{R}^2$, $\lambda,\delta,{\delta'}\in (0,1]$ and $\phi \in \mathcal{B}_0$.
\end{Lemma}
\begin{proof}
Unraveling \eqref{eq:PAM prod guys}, we find that
\begin{align*}
\hat{\Pi}_\star^\delta \partial_i I \Xi \cdot \partial_j I\Xi \bigl(\phi^\lambda_\star\bigr)&{}= {\Pi}_\star^\delta \partial_i I \Xi \cdot \partial_j I\Xi \bigl(\phi^\lambda_\star\bigr)
- c^{\<b2>_{ij}}_\delta \bigl(\phi^\lambda_\star\bigr) \\
&{}= \int_{(\mathbb{R}^2)^3} \partial_i H(x,y_1) \partial_i H(x,y_2) \xi_\delta(y_1)\xi_\delta(y_2) \phi^\lambda_\star(x) \,\dd y_1 \dd y_2 \dd x
- c^{\<b2>_{ij}}_\delta \bigl(\phi^\lambda_\star\bigr).
\end{align*}
As in the proof of Lemma~\ref{lem:dumbbell}, for $q\in \mathbb{N}$, let $\psi\colon \mathbb{R}^2\to [0,1]$ be smooth, compactly supported on ${B}_{1/2q}(0)\subset \mathbb{R}^2$ such that
\smash{$\psi_q|_{{B}_{1/4q}(0)}=1$}. This time set
$K_i(x,y)=\psi(x-y) \partial_i H(x,y)$,
\smash{$R_i(x,y):= \partial_i H(x,y)-K_i(x,y)$} and
\[
\hat{c}_{\delta}^{\<b2>_{ij}}(x) = \int_{(\mathbb{R}^2)^2} \psi(y_1) \partial_i H(x,y_1) \psi(y_1) \partial_j H(x,y_2) \rho^\delta *\rho^\delta (y_1-y_2)\, \dd y_1\dd y_2.
\]
Thus,
\begin{align*}
\hat{\Pi}_\star^\delta \partial_i I \Xi \cdot \partial_j I\Xi \bigl(\phi^\lambda_\star\bigr)
 ={}& \int_{(\mathbb{R}^2)^3} \partial_i H(x,y_1) \partial_i H(x,y_2) \xi_\delta(y_1)\xi_\delta(y_2) \phi^\lambda_\star(x)\, \dd y_1 \dd y_2 \dd x\\
 &{}
{-}\, \hat{c}^{\<b2>_{ij}}_\delta \bigl(\phi^\lambda_\star\bigr) -\bigl(c^{\<b2>_{ij}}_\delta -\hat{c}^{\<b2>_{ij}}_\delta\bigr)\bigl(\phi^\lambda_\star\bigr) \\
 ={}&\int_{(\mathbb{R}^2)^3} K_i (x,y_1) K_j (x,y_1) \xi_\delta(y_1)\xi_\delta(y_2) \phi^\lambda_\star(x)\, \dd y_1 \dd y_2 \dd x
 - \hat{c}^{\<b2>_{ij}}_\delta \bigl(\phi^\lambda_\star\bigr)
 \\
 &{}{+}\,
\int_{(\mathbb{R}^2)^3} R_i (x,y_1) K_j (x,y_1) \xi_\delta(y_1)\xi_\delta(y_2) \phi^\lambda_\star(x)\, \dd y_1\dd y_2 \dd x
 \\
 &{}{+}\,
\int_{(\mathbb{R}^2)^3} K_i (x,y_1) R_j (x,y_1) \xi_\delta(y_1)\xi_\delta(y_2) \phi^\lambda_\star(x)\, \dd y_1 \dd y_2 \dd x
 \\
 &{}{+}\,
\int_{(\mathbb{R}^2)^3} R_i(x,y_1) R_j (x,y_1) \xi_\delta(y_1)\xi_\delta(y_2) \phi^\lambda_\star(x)\, \dd y_1 \dd y_2 \dd x\\
&{}{-}\,\bigl(c^{\<b2>_{ij}}_\delta -\hat{c}^{\<b2>_{ij}}_\delta\bigr)\bigl(\phi^\lambda_\star\bigr).
\end{align*}
Estimating the first term is the content of Lemma~\ref{lem:cherry}. The two last terms are easily bounded directly.
Finally, for the remaining terms one notes as a consequence of \cite[Lemma~2.9]{Sin25}
that
$R_i\in \bfK^{2}_{L-1,R}$ almost surely and that the random variables $\| R_i\|_{\beta; L-1,R}$ have finite $q$-th moment for any $q<\infty$.
Then, as in the proof of Lemma~\ref{lem:dumbbell}, it follows that
\[ \sup_{\delta\in [0,1]}\biggl\| \int_{\mathbb{R}^2} R_i(\cdot,y) \xi_\delta(y)\biggr\|_{C^{L-1}} \in L^{q}(\Omega)\] and the desired bound follows by Young multiplication.

The difference bound for two values $\delta,\delta'\in (0,1]$ follows similarly.
\end{proof}

 \begin{Lemma}\label{lem:cherry}
 For $\alpha<0$, $q\in 2\mathbb{N}$ and $\kappa\in (0,1)$ in the setting of the proof of Lemma~{\rm \ref{lem:dumbbell2}}, it holds that for $n,m\in\{1, 2\}$
\[\E \biggl[ \biggl|\int_{(\mathbb{R}^2)^3} K_n(x,y) K_m(x,y') \bigl[\xi_\delta(y)\diamond \xi_\delta(y')\bigr] \phi^\lambda_\star(x)\,\dd y \dd y' \dd x \biggl|^q \biggr]^{\frac{1}{q}}\lesssim_{q,\alpha} \lambda^\alpha,\]
as well as
\begin{align*}
&\E \biggl[ \biggl|\int_{(\mathbb{R}^2)^3} K_n(x,y) K_m(x,y') \bigl[\xi_{\delta}(x)\diamond\xi_{\delta}(y)-\xi_{{\delta'}}(x)\diamond \xi_{{\delta'}}(y)\bigr] \phi^\lambda_\star(x)\,\dd y \dd y' \dd x \biggr|^q \biggr]^{\frac{1}{q}}\\
&\qquad{}
\lesssim_{q,\alpha} |\delta-{\delta'}|^\kappa \lambda^{\alpha-\kappa}
\end{align*}
uniformly over $\star\in \mathbb{R}^2$, $\lambda,\delta, \delta' \in (0,1]$ and $\phi \in \mathcal{B}_0$.

 \end{Lemma}
 \begin{proof}
As previously, to exhibit the dependence on the Gaussian explicitly, set for $n\in \{1,2\}$
\[\tilde{K}^\eta_n(x,y)=\psi(x-y) \tilde{G}_n^{\eta}(x-y).\]
 Then,
by Gaussian integration by parts, see Lemma~\ref{lem:wick_integration by parts}, one finds that
 \begin{align*}
&\E \biggl[ \biggr|\int_{(\mathbb{R}^2)^3} K_n(x,y) K_m(x,y') \bigl[\xi_\delta(y)\xi_\delta(y')- \E\bigl[\xi_\delta(y)\xi_\delta(y')\bigr]\bigr] \phi^\lambda_\star(x)\,\dd y \dd y' \dd x \biggr|^q \biggr]\\
&\qquad =\E \biggl[ \biggl|\int_{(\mathbb{R}^2)^3}\tilde{K}^{h(x)}_n (x,y) \tilde{K}^{h(x)}_m(x,y') \bigl[\xi_\delta(y)\xi_\delta(y')- \E\bigl[\xi_\delta(y)\xi_\delta(y')\bigr]\bigr] \phi^\lambda_\star(x)\,\dd y \dd y' \dd x \biggr|^q \biggr]\\
&\qquad =
\sum_{
{\underset{\vert J\vert\text{ even}}
{J\subset\{1,\dots,2q\}}}}
\sum_{P\in \mathcal{P}_2(J)}
\int\Biggl[
\prod_{\{i,j\}\in P}\mathbb{E}[\xi_\delta(y_i)\xi_\delta(y_j)]
\prod_{i=1}^q \phi^\lambda_\star(x_i)
 \\
&\hphantom{\qquad =\sum_{{\underset{\vert J\vert\text{ even}}{J\subset\{1,\dots,2q\}}}}\sum_{P\in \mathcal{P}_2(J)}\int\Biggl[}{}
\times
\sum_{k\colon J^c\rightarrow \{1,\dots,q\}}\!\!
\mathbb{E}\Biggl[\prod_{i=1}^{q} \partial_{\uparrow}^{|k^{-1}(i)|} \bigl(\tilde{K}^{h(x_i)}_n(x_i,y_{2i-1}) \tilde{K}^{h(x_i)}_m(x_i,y_{2i})\bigr)\Biggr]\\
&\hphantom{\qquad =\sum_{{\underset{\vert J\vert\text{ even}}{J\subset\{1,\dots,2q\}}}}\sum_{P\in \mathcal{P}_2(J)}\int\Biggl[}{}
\times
\prod_{i\in J^c}\mathbb{E}\bigl[\xi_\delta(y_i) \xi*\sigma \bigl({x_{k(i)}}\bigr) \bigr] \Biggr],
\end{align*}
where the integral is over \smash{$(y_1,\dots,y_{2q})\in\bigl(\mathbb{R}^{2}\bigr)^{2q}$} and \smash{$(x_1,\dots,x_q)\in \bigl(\mathbb{R}^{2}\bigr)^{q}$}.
Inserting \eqref{eq:correlation identities}, this can be bounded by
\begin{align*}
&\sum_{J, P, k, \ell}\Biggl|
\int
\prod_{\{i,j\}\in P}
\bigl(\rho^{*2}\bigr)^\delta(y_i-y_j)
 \prod_{i=1}^q \phi^\lambda_\star(x_i)\\
& \hphantom{\sum_{J, P, k, \ell}\Biggl|\int}{}
\times
\mathbb{E}\Biggl[\prod_{i=1}^{q} \partial_{\uparrow}^{|k^{-1}(i)|-\ell_i} \tilde{K}^{h(x_i)}_n(x_i,y_{2i-1}) \partial_{\uparrow}^{\ell_i} \tilde{K}^{h(x_i)}_m(x_i,y_{2i})\Biggr]\prod_{i\in J^c}
F^{x_{k(i)}}_\delta(y_i)\Biggr|,
\end{align*}
where the sum runs over even subsets $J\subset\{1,\dots,2q\}$, partitions $P\in \mathcal{P}_2(J)$, maps
$k \colon J^c\rightarrow \{1,\dots,q\}$
and $\ell\in \mathbb{N}^q$ such that
\smash{$|\ell_i|<\bigl|k^{-1}(i)\bigr|$} for each $i\in \{1,\dots,q\}$.
Thus, we shall show that for each summand
\begin{gather*}
\mathbb{E}\Biggl[\Biggl|
\int
\bigl(\rho^{*2}\bigr)^\delta(y_i-y_j)
 \prod_{i=1}^q \phi^\lambda_\star(x_i)
\prod_{i=1}^{q} \partial_{\uparrow}^{|k^{-1}(i)|-\ell_i} \tilde{K}^{h(x_i)}_n(x_i,y_{2i-1}) \partial_{\uparrow}^{\ell_i} \tilde{K}^{h(x_i)}_m(x_i,y_{2i})\\
\hphantom{\mathbb{E}\Biggl[\Biggl|\int}{}
\times\prod_{i\in J^c}
F^{x_{k(i)}}_\delta(y_i)\Biggr|\Biggr] \lesssim \lambda^{q \alpha},
\end{gather*}
using the modified Hairer--Quastel criterion, Theorem~\ref{thm:slight variant of HQ}.
We represent the integral using a~directed graph $\CCG=(\CCV,\CCE)$ constructed as follows.
The vertex set is given by $\CCV =\{ v_0,\dots,v_{q},\allowbreak w_1,\dots, w_{2q}\}$. Each element $v_i$ for $i\geq 1$ of which is interpreted as representing the integration variable~$x_i$,~$v_0$ as representing $\star$ and $w_j$ as representing $y_j$.
 \begin{itemize}\itemsep=0pt
\item It contains $q$ distinguished edges representing the generic test function $\phi_\star^\lambda$, which we shall illustratively draw as
 \tikz[baseline=-0.1cm] \draw[testfcn] (1,0) to (0.0,0) \enlarge;, one directed from $v_0$ to each odd vertex $v_{i}$.
 We set ${(a_e,r_e)=(0, 0)}$.
\item It contains $q$ edges representing factors of
\smash{$\partial_{\uparrow}^{|k^{-1}(i)|-\ell_i} \tilde{K}^{h(x_i)}_n(x_i,y_{2i-1}) $}
for each $i=1,\dots,q$, one directed the each odd vertex $w_{2i-1}$ to $v_{i}$ as well as another $q$ edges representing a~factors of
\smash{$\partial_{\uparrow}^{\ell_i} \tilde{K}^{h(x_i)}_m(x_i,y_{2i})$} for each $i=1,\dots,q$, one directed from the even vertex $w_{2i}$ to~$v_{i}$.
We set $(a_e, r_e)=(1+\kappa',0)$ for some $\kappa'\in (0,1)$ specified later. We call these kernel edges and draw then as \tikz[baseline=-0.1cm] \draw[kernel] (0.0,0) to (1,0);.
\item It contains $|J|/2$ edges representing $\bigl(\rho^{*2}\bigr)^\delta$, drawn as \tikz[baseline=-0.1cm] \draw[rho] (0.0,0) to (1,0); (we again suppress the arbitrary orientation since it is irrelevant). For each $P\in \mathcal{P}_2(J)$ and $\{i,j\} \in P$ it connects $v_i$ to $v_j$.
 We set $I_e=1$ for these edges, as well as $(a_e,r_e)=(2+\kappa, -1)$ for some $\kappa\in (0,1)$ specified later and call these $\rho$-type edges.
 \item Finally, there are $|J^c|$ edges, for each $i\in J^c$ one connecting $v_{k(i)}$ to $w_i$ representing the function \smash{$F^{x_{k(i)}}_\delta(x_i)$} for which we set $(a_e, r_e)=(0,0)$. (We again omit them from the schematic drawings below).
 \end{itemize}

 We draw some examples (up to the omitted data) for $q=2$ and $q=4$
 
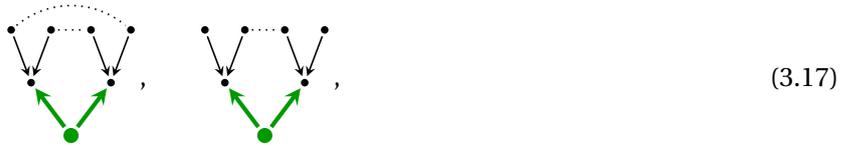
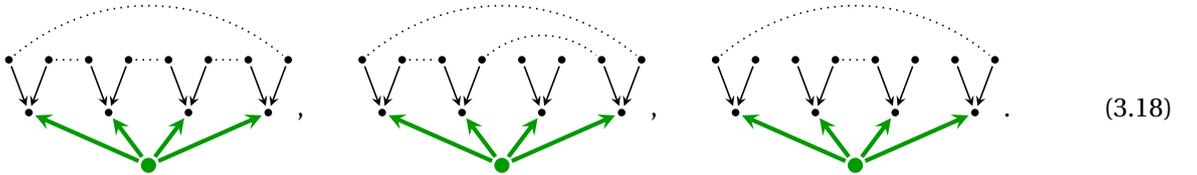
\begin{gather*}
 \begin{tikzpicture}[scale=0.35,baseline=0.3cm]
	\node at (0,-1) [root] (root) {};
	\node at (-1.5,1) [dot] (a) {};
	\node at (1.5,1) [dot] (b) {};
	\node at (-2.25,3) [dot] (aa) {};
	\node at (-0.75,3) [dot] (ab) {};
	\node at (2.25,3) [dot] (bb) {};
	\node at (0.75,3) [dot] (ba) {};
	\draw[testfcn] (a) to (root);
	\draw[testfcn] (b) to (root);
	\draw[kernel] (aa) to (a);
	\draw[kernel] (ab) to (a);
	\draw[kernel] (ba) to (b);
	\draw[kernel] (bb) to (b);
	\draw[rho] (aa) to[bend left=40] (bb);
	\draw[rho] (ab) to (ba);
\end{tikzpicture}, \qquad
 \begin{tikzpicture}[scale=0.35,baseline=0.3cm]
	\node at (0,-1) [root] (root) {};
	\node at (-1.5,1) [dot] (a) {};
	\node at (1.5,1) [dot] (b) {};
	\node at (-2.25,3) [dot] (aa) {};
	\node at (-0.75,3) [dot] (ab) {};
	\node at (2.25,3) [dot] (bb) {};
	\node at (0.75,3) [dot] (ba) {};
	\draw[testfcn] (a) to (root);
	\draw[testfcn] (b) to (root);
	\draw[kernel] (aa) to (a);
	\draw[kernel] (ab) to (a);
	\draw[kernel] (ba) to (b);
	\draw[kernel] (bb) to (b);
	\draw[rho] (ab) to (ba);
\end{tikzpicture},
\\ 
\begin{tikzpicture}[scale=0.35,baseline=0.3cm]
	\node at (0,-1) [root] (root) {};
	\node at (-4.5,1) [dot] (a) {};
	\node at (-1.5,1) [dot] (b) {};
	\node at (1.5,1) [dot] (c) {};
	\node at (4.5,1) [dot] (d) {};	
	\node at (-5.25,3) [dot] (aa) {};
	\node at (-3.75,3) [dot] (ab) {};	
	\node at (-2.25,3) [dot] (ba) {};
	\node at (-0.75,3) [dot] (bb) {};
	\node at (2.25,3) [dot] (cb) {};
	\node at (0.75,3) [dot] (ca) {};
	\node at (5.25,3) [dot] (db) {};
	\node at (3.75,3) [dot] (da) {};	
	\draw[testfcn] (a) to (root);
	\draw[testfcn] (b) to (root);
	\draw[testfcn] (c) to (root);
	\draw[testfcn] (d) to (root);
	\draw[kernel] (aa) to (a);
	\draw[kernel] (ab) to (a);
	\draw[kernel] (ba) to (b);
	\draw[kernel] (bb) to (b);
	\draw[kernel] (ca) to (c);
	\draw[kernel] (cb) to (c);
	\draw[kernel] (da) to (d);
	\draw[kernel] (db) to (d);
	\draw[rho] (aa) to[bend left=40] (db);
	\draw[rho] (ab) to (ba);
	\draw[rho] (bb) to (ca);
	\draw[rho] (cb) to (da);
\end{tikzpicture}, \qquad
\begin{tikzpicture}[scale=0.35,baseline=0.3cm]
	\node at (0,-1) [root] (root) {};
	\node at (-4.5,1) [dot] (a) {};
	\node at (-1.5,1) [dot] (b) {};
	\node at (1.5,1) [dot] (c) {};
	\node at (4.5,1) [dot] (d) {};	
	\node at (-5.25,3) [dot] (aa) {};
	\node at (-3.75,3) [dot] (ab) {};	
	\node at (-2.25,3) [dot] (ba) {};
	\node at (-0.75,3) [dot] (bb) {};
	\node at (2.25,3) [dot] (cb) {};
	\node at (0.75,3) [dot] (ca) {};
	\node at (5.25,3) [dot] (db) {};
	\node at (3.75,3) [dot] (da) {};	
	\draw[testfcn] (a) to (root);
	\draw[testfcn] (b) to (root);
	\draw[testfcn] (c) to (root);
	\draw[testfcn] (d) to (root);
	\draw[kernel] (aa) to (a);
	\draw[kernel] (ab) to (a);
	\draw[kernel] (ba) to (b);
	\draw[kernel] (bb) to (b);
	\draw[kernel] (ca) to (c);
	\draw[kernel] (cb) to (c);
	\draw[kernel] (da) to (d);
	\draw[kernel] (db) to (d);
	\draw[rho] (aa) to[bend left=40] (db);
	\draw[rho] (ab) to (ba);
	\draw[rho] (bb) to[bend left=40] (da);
\end{tikzpicture}, \qquad
\begin{tikzpicture}[scale=0.35,baseline=0.3cm]
	\node at (0,-1) [root] (root) {};
	\node at (-4.5,1) [dot] (a) {};
	\node at (-1.5,1) [dot] (b) {};
	\node at (1.5,1) [dot] (c) {};
	\node at (4.5,1) [dot] (d) {};	
	\node at (-5.25,3) [dot] (aa) {};
	\node at (-3.75,3) [dot] (ab) {};	
	\node at (-2.25,3) [dot] (ba) {};
	\node at (-0.75,3) [dot] (bb) {};
	\node at (2.25,3) [dot] (cb) {};
	\node at (0.75,3) [dot] (ca) {};
	\node at (5.25,3) [dot] (db) {};
	\node at (3.75,3) [dot] (da) {};	
	\draw[testfcn] (a) to (root);
	\draw[testfcn] (b) to (root);
	\draw[testfcn] (c) to (root);
	\draw[testfcn] (d) to (root);
	\draw[kernel] (aa) to (a);
	\draw[kernel] (ab) to (a);
	\draw[kernel] (ba) to (b);
	\draw[kernel] (bb) to (b);
	\draw[kernel] (ca) to (c);
	\draw[kernel] (cb) to (c);
	\draw[kernel] (da) to (d);
	\draw[kernel] (db) to (d);
	\draw[rho] (aa) to[bend left=40] (db);
	\draw[rho] (bb) to (ca);
\end{tikzpicture}.
\end{gather*}

Now we check the assumptions of the modified Hairer--Quastel criterion, Theorem~\ref{thm:slight variant of HQ}.
The first item of Assumption~\ref{ass:mainGraph_weak}, i.e., that for every edge $e\in \CCE$, one has $a_{e}-r^{-}_e<2$ is clearly satisfied.
The second item, i.e., \eqref{e:assEdges}, is easily checked by a similar case separation as in the proof of Lemma~\ref{lem:old_term}.
Thus, it remains to calculate the quantities
\[
\tilde \alpha = |\fraks||\CCV\setminus \CCV_\star| - \sum_{e\in \CCE} a_e-R(\CCG), \qquad
 R(\CCG):= \max_{\bar \CCV \subset \CCV\setminus \CCV_\star} \biggl( |\bar \CCV||\fraks|- \sum_{e\in \CCE(\bar\CCV) } a_e\biggr)\vee 0.
 \]
Let $m$ be the number of vertices in $\CCV\setminus \CCV_\star$ which have only one incident edge. Then one sees choosing $\bar \CCV$ to be this subset maximises the
$R(\CCG)= 2m - m(1+\kappa')= m(1-\kappa')$.
Since \smash{$\sum_{e\in \CCE} a_e= 2q(1+\kappa') + (q-m/2) (2+\kappa) $}, we find that
\begin{align*}
 \tilde \alpha = |\fraks||\CCV\setminus \CCV_\star| - \sum_{e\in \CCE} a_e-R(\CCG)
&{}= 4q- 2q(1+\kappa') - (q-m/2) (2+\kappa) -m(1-\kappa)\\
&{}=- q(2 \kappa' + \kappa) +3m\kappa/2.
\tag*{\qed}
\end{align*}
\renewcommand{\qed}{}
\end{proof}

\begin{Remark}
It might be suggestive to think that one can simply integrate out the substructures%
 
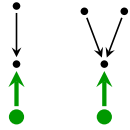
\begin{equation*}
 \begin{tikzpicture}[scale=0.35,baseline=0.3cm]
	\node at (0,-1) [root] (root) {};
	\node at (0,1) [dot] (b) {};
	
	\node at (0,3.2) [dot] (bb) {};

	\draw[testfcn] (b) to (root);
	\draw[kernel] (bb) to (b);
\end{tikzpicture}, \qquad
 \begin{tikzpicture}[scale=0.35,baseline=0.3cm]
	\node at (0,-1) [root] (root) {};
	\node at (0,1) [dot] (b) {};
	
	\node at (0.75,3) [dot] (bb) {};
	\node at (-0.75,3) [dot] (ba) {};

	\draw[testfcn] (b) to (root);
	
	\draw[kernel] (ba) to (b);
	\draw[kernel] (bb) to (b);
\end{tikzpicture}
 \end{equation*}
 whenever they appear, but this is not obvious due to the functions $F$ hidden in the graphical schematic.
\end{Remark}

\section[Stochastic estimates for the phi\^{}{K+1}\_2-equation]{Stochastic estimates for the $\boldsymbol{\phi^{K+1}_2}$-equation}
In this section, we prove the following lemma.
\begin{Lemma}\label{lem:phi4 estimate}
{\samepage For any $\alpha<0$, $N\in\mathbb{N}$ and $q\in 2\mathbb{N}$, there exists $\kappa>0$ such that
\begin{gather*}
\E\bigl[ \bigl|H_N\bigl(\<1>_\delta, c^{\<2>}_\delta \bigr)\bigl(\phi^\lambda_\star\bigr)\bigr|^q\bigr]\lesssim_{q,N} \lambda^{q \alpha}
\qquad \text{and} \\
\E\bigl[
\bigl|H_N\bigl(\<1>_\delta, c^{\<2>}_\delta \bigr)\bigl(\phi^\lambda_\star\bigr)- H_N\bigl(\<1>_{\delta'}, c_{\delta'}^{\<2>} \bigr)\bigl(\phi^\lambda_\star\bigr) \bigr|^q
\bigr]\lesssim_{\alpha,q} |\delta-{\delta'}|^{\kappa q} \lambda^{(\alpha -\kappa)q},
\end{gather*}
 uniformly over $\star\in \mathbb{R}^2$, $\lambda,\delta,{\delta'}\in (0,1]$ and $\phi \in \mathcal{B}_0$.}
\end{Lemma}
In the proof, we shall use the following identity about Hermite polynomials, which is checked using the recursive definition \eqref{eq:Hermit}
\begin{gather}\label{eq:Hermit identity}
H_{N}(X+Y,c+d)= \sum_{n=0}^N {N \choose n } H_{N-n}(X,c) H_{N}(Y,d).
\end{gather}
\begin{proof}[Proof of Lemma~\ref{lem:phi4 estimate}]
Similarly to the proofs in Section~\ref{sec:stoch est pam}, let $\psi\colon \mathbb{R}^2\to [0,1]$ be smooth, compactly supported on ${B}_{1/4Nq}(0)\subset \mathbb{R}^2\times $ such that
$\psi\big|_{{B}_{1/8Nq}(0)}=1$ and set
\begin{gather*}
K(t,x;s,y)=\psi(x-y) \kappa(t-s) Z^{(t,x)}(t-s,x-y),\\
R(t,x;s,y):= \kappa(t-s) \Gamma(t,x;s,y)-K(t,x;s,y)
\end{gather*}
and \[
\hat{\<1>}_\delta(z)= \int K(z,w) \xi_\delta(w)\, \dd w, \qquad
\hat{c}^{\<2>}_{\delta}(z) = \int \int K(z;w)K(z;w') \bigl(\rho^\delta *\rho^\delta\bigr) (w-w')\, \dd w \dd w'. \]
Thus by \eqref{eq:Hermit identity}
\[ H_N\bigl(\<1>_\delta, c_\delta^{\<2>} \bigr)\bigl(\phi^\lambda_\star\bigr)=
 \sum_{n=0}^N {N \choose n } H_{N-n}\bigl(\hat{\<1>}_\delta,\hat{c}^{\<2>}_\delta\bigr) H_{n} \bigl({\<1>}_\delta-\hat{\<1>}_\delta,c^{\<2>}_\delta-\hat{c}^{\<2>}_\delta\bigr).
\]
Noting that \smash{${\<1>}_\delta-\hat{\<1>}_\delta= \int R(z,w) \xi_\delta(w)\, \dd w$}, it follows from Levi's construction of the fundamental solution, cf.\ \cite[Proposition~2.6]{Sin25} that
\smash{$R\in \bfK^{3}_{L,R}$} for $R\in (2,3)$ and $L>0$ almost surely and that the random variable $\|R\|_{3; L,R}$ has a finite $q$-th moment for any $q<\infty$. Since the variable $ \sup_{\delta\in [0,1]} \| \xi_\delta(x)\|_{C^{-2^{-}}}$ also has finite $q$-th moment for any $q<\infty$, it follows that
there exists $\kappa>0$ such that for any $q<\infty$
\begin{gather}\label{eq:higher regularity used}
 \sup_{\delta} \bigl\|{\<1>}_\delta-\hat{\<1>}_\delta\bigr\|_{C^{\kappa}} \in L^{q}(\Omega).
 \end{gather}
Combining this with classical Young multiplication and Lemma~\ref{lem:phi4 singular estimate} completes the proof.
The difference bound follows similarly.
\end{proof}

\begin{Lemma}\label{lem:phi4 singular estimate}
For $\alpha<0$, $q\in 2\mathbb{N}$ and $\kappa\in (0,1)$, in the setting of the proof of Lemma~{\rm \ref{lem:phi4 estimate}}, it holds that
\begin{gather*}
\E\bigl[ \bigl|H\bigl(\hat{\<1>}_\delta, \hat{c}^{\<2>}_\delta \bigr)\bigl(\phi^\lambda_\star\bigr)\bigr|^q\bigr]\lesssim_{q,N} \lambda^{q \alpha}
\qquad \text{and}\\
\E\bigl[ \bigl|H\bigl(\hat{\<1>}_\delta, \hat{c}^{\<2>}_\delta \bigr)\bigl(\phi^\lambda_\star\bigr)- H\bigl(\hat{\<1>}_{\delta'}, \hat{c}^{\<2>}_{\delta'} \bigr)\bigl(\phi^\lambda_\star\bigr)\bigr|^q \bigr]\lesssim_{q,N} |\delta- \delta'|^{q\kappa} \lambda^{q(\alpha-\kappa)}
\end{gather*}
 uniformly over $\star\in \mathbb{R}^2$, $\lambda,\delta, \delta' \in (0,1]$ and $\phi \in \mathcal{B}_0$.
\end{Lemma}
\begin{proof}
Recalling the notion of Wick products, cf.\ Appendix~\ref{App:Gaussian int}, we first check the following identity for $N\in \mathbb{N}$:
\begin{gather}\label{WickproductHermit}
H_N\bigl(\hat{\<1>}(x), \hat{c}^{\<2>}_\delta(x)\bigr)=\int \prod_{i=1}^N K(x,y_i)\biggl[\mathop{\diamond}_{i=1}^N \xi_\delta(y_i)\biggr]\qquad\text{for $x\in \mathbb{R}^2$}
\end{gather}
by induction, where the integral is over \smash{$(y_1,\dots,y_N)\in \bigl(\mathbb{R}^3\bigr)^{N}$}.
Indeed, for $N\in \{0,1\}$, \eqref{WickproductHermit}~holds. For $N\geq 2$, using the definition of $H_N$ and the induction hypothesis, one finds that
\begin{align*}
\begin{split}
H_N\bigl(\<1>(x), c^{\<2>}_\delta(x)\bigr)={}& \<1>(x)H_{N-1}\bigl(\<1>(x), c^{\<2>}_\delta(x)\bigr)-(N-1)c_\delta(x)H_{N-2}\bigl(\<1>(x), c^{\<2>}_\delta(x)\bigr)\\
={}&\int \prod_{i=1}^N K(x,y_i) \xi_\delta(y_N)\mathop{\diamond}_{i=1}^{N-1} \xi_\delta(y_i)\\
&{}{-}\,(N-1)\int \prod_{i=1}^{N}K(x,y_i) \mathbb{E}[\xi_\delta(y_{N-1})\xi_\delta(y_N)]\mathop{\diamond}_{i=1}^{N-2}\xi_\delta(y_i).
\end{split}
\end{align*}
Using the symmetry \eqref{eq:permutation} of the Wick product, we note that
\begin{align*}
&(N-1)\int \prod_{i=1}^{N}K(x,y_i)\, \mathbb{E}[\xi_\delta(y_{N-1})\xi_\delta(y_N)]\mathop{\diamond}_{i=1}^{N-2}\xi_\delta(y_i)\\
&\qquad =\int \prod_{i=1}^{N}K(x,y_i)\, \sum_{\ell=1}^{N-1}\mathbb{E}[\xi_\delta(y_{\ell})\xi_\delta(y_N)]\mathop{\diamond}_{\underset{j\neq \ell}{j=1}}^{N-1}\xi_\delta(y_i)
\end{align*}
and thus obtain
\[
H_{N}(\<1>(x), c_\delta(x))=\int \prod_{i=1}^N K(x,y_i) \Biggl(\xi_\delta(y_N)\mathop{\diamond}_{i=1}^{N-1} \xi_\delta(y_i)-\sum_{\ell=1}^{N-1}\mathbb{E}[\xi_\delta(y_{\ell})\xi_\delta(y_N)]\mathop{\diamond}_{\underset{j\neq \ell}{j=1}}^{N-1}\xi_\delta(y_i)\Biggr)
\]
and therefore \eqref{WickproductHermit}.
To exhibit the dependence of the kernel on the Gaussian explicitly, set
\[\tilde{K}^\eta(t,\zeta;s,\chi)=\psi(\zeta-\chi) \kappa(t-s) \tilde{Z}^{\eta}(t-s,\zeta-\chi)\;
\]
for $\tilde{Z}$ defined in \eqref{eq:frozen kernel_factored}.
Thus, for $q\in 2\mathbb{N}$
\begin{align*}
&\E\bigl[ \bigl|H\bigl(\hat{\<1>}_\delta, \hat{c}^{\<2>}_\delta \bigr)\bigl(\phi^\lambda_\star\bigr)\bigr|^q\bigr]\\
&\qquad=\E\Biggl[ \Biggl|\int \phi^\lambda_\star (x)\prod_{i=1}^N \tilde{K}^{h(x)}(x, z_i)\biggl[\mathop{\diamond}_{i=1}^N \xi_\delta(z_i) \biggr]\Biggr|^q \Biggr]\\
&\qquad=
\int \E\Biggl[
\prod_{j=1}^{q}\Biggl[
\phi^\lambda_\star(x_j)
\prod_{i=1}^N
\tilde{K}^{h(x_j)}(x_j, z_{(j-1)N+i}) \Biggr]
\prod_{j=1}^{q} \biggl[\mathop{\diamond}_{i=1}^N \xi_\delta(z_{(j-1)N+i}) \biggr] \Biggr],
\end{align*}
where the integral in the last term is over \smash{$(z_1,\dots,z_{qN})\in \bigl(\mathbb{R}^3\bigr)^{qN}$} and \smash{$(x_1,\dots,x_N)\in \bigl(\mathbb{R}^3\bigr)^{N}$}.
Next, we use Proposition~\ref{prop:wickInteByPartsGeneral}
to find that this agrees with
\begin{align*}
&\sum_{
{\underset{\vert J\vert\text{ even}}
{J\subset\{1,\dots,Nq\}}}
}\sum_{P\in \mathcal{P}_N(J)} \int \Biggl[
\prod_{j=1}^{q}
\phi^\lambda_\star(x_j)
\prod_{\{i,j\}\in P}
\mathbb{E}[\xi_\delta(z_i)\xi_\delta(z_j)]
\\
&\hphantom{\sum_{{\underset{\vert J\vert\text{ even}}{J\subset\{1,\dots,Nq\}}}}\sum_{P\in \mathcal{P}_N(J)} \int\Biggl[}{}
\times\sum_{k \colon J^c\rightarrow \{1,\dots,q\}}
\mathbb{E}\Biggl[\prod_{j=1}^q \partial_{\uparrow}^{|k^{-1}(j)|} \prod_{i=1}^N
\tilde{K}^{h(x_j)}\bigl(x_j, z_{(j-1)N+i}\bigr)\Biggr]\\
&\hphantom{\sum_{{\underset{\vert J\vert\text{ even}}{J\subset\{1,\dots,Nq\}}}}\sum_{P\in \mathcal{P}_N(J)} \int\Biggl[}{}
\times\prod_{i\in J^c}\mathbb{E}\bigl[\xi_{\delta}(z_i) h\bigl(x_{k(i)}\bigr)\bigr] \Biggr],
\end{align*}
where we recall that
\begin{gather*}
\mathcal{P}_N(J)=\{ P\in \mathcal{P}(J)\mid |\{Ni\!+\!1,Ni\!+\!2,\dots, N(i\!+\!1) \} \cap p| \leq 1 \, \forall i\in \{0,\dots,q-1\}, \, \forall p \in P \}.
\end{gather*}
Using the obvious space-time analogues of \eqref{eq:correlation identities},
this in turn equals to
\begin{align*}
&\sum_{J,P,k}
\int
\prod_{j=1}^{q}
\phi^\lambda_\star(x_j)
\prod_{\{i,j\}\in P}
\bigl(\rho^{*2}\bigr)^\delta(z_i-z_j)\\
&\hphantom{\sum_{J,P,k}\int}{}\times
\mathbb{E}\Biggl[\prod_{j=1}^q \partial_{\uparrow}^{|k^{-1}(j)|}\prod_{i=1}^N
\tilde{K}^{h(x_j)}\bigl(x_j, z_{(j-1)N+i}\bigr)\Biggr]\prod_{i\in J^c}F_\delta^{x_{k(i)}}(z_i)\\
&\qquad{}\lesssim
\sum_{J,P,k}
\sum_{\substack{
\ell^i\in \mathbb{N}^N\ :\\ \ |\ell^i|=|k^{-1}(j)|} }
\Biggl|\int
\prod_{j=1}^{q}
\phi^\lambda_\star(x_j)
\prod_{\{i,j\}\in P}
\bigl(\rho^{*2}\bigr)^\delta(z_i-z_j)\\
&\qquad\hphantom{\lesssim\sum_{J,P,k}\sum_{\substack{\ell^i\in \mathbb{N}^N\ :\\ \ |\ell^i|=|k^{-1}(j)|} }\Biggl|\int}{} \times
\mathbb{E}\Biggl[\prod_{j=1}^q\prod_{i=1}^N
\partial_{\uparrow}^{\ell_{i}^j}
\tilde{K}^{h(x_j)}\bigl(x_j, z_{(j-1)N+i}\bigr)\Biggr]\prod_{i\in J^c}F_\delta^{x_{k(i)}}(z_i)\Biggr|,
\end{align*}
where the sums over $J$, $P$, $k$ are over even subsets $J\subset \{1,\dots,Nq\}$, pairings $P\in \mathcal{P}_N(J)$ and maps $ J^c\rightarrow \{1,\dots,q\} $.
Thus, we shall bound each such summand
\begin{gather*}
\E\Biggl[\Biggl|\int
\prod_{j=1}^{q}
\phi^\lambda_\star(x_j)
\prod_{\{i,j\}\in P}
\bigl(\rho^{*2}\bigr)^\delta(z_i-z_j)\\
\hphantom{\E\Biggl[\Biggl|\int}{}\times
\sum_{k\colon J^c\rightarrow \{1,\dots,q\}}
\prod_{i=1}^N
\partial_{\uparrow}^{\ell_{i}^j}
\tilde{K}^{h(x_j)}\bigl(x_j, z_{(j-1)N+i}\bigr)\prod_{i\in J^c}F_\delta^{x_{k(i)}}(z_i) \Biggr|\Biggr]
\end{gather*}
using a directed graph $\CCG=(\CCV,\CCE)$ constructed as follows.
The vertex set is given by $\CCV =\{ v_0,\dots,v_{q}, w_1,\dots w_{Nq}\}$, each element $v_i$ for $i\geq 1$ of which we interpret as representing the integration variable $x_i$ and $v_0$ representing $\star$ as well as $w_j$ representing $z_j$.
 \begin{itemize}\itemsep=0pt
\item It contains $q$ distinguished edges representing the generic test function $\phi_\star^\lambda$, one directed from $v_0$ to each vertex $v_{i}$.
 We set $(a_e,r_e)=(0, 0)$ and illustratively draw these edges as~\tikz[baseline=-0.1cm] \draw[testfcn] (1,0) to (0.0,0) \enlarge;.
\item It contains $Nq$ $K$-type edges, one for each $(i,j)\in\{1,\dots,q\}\times \{1,\dots,N\}$ representing a~factors
\smash{$\tilde{K}^{h(x_j)}\bigl(x_j, z_{(j-1)N+i}\bigr)$}
which is directed from $w_{(j-1)N+i}$ to $v_{j}$.
We all draw these as~\tikz[baseline=-0.1cm] \draw[kernel] (0.0,0) to (1,0);. We set $(a_e, r_e)=(2,0)$.
\item It contains $|J|/2$ edges representing \smash{$\bigl(\rho^{*2}\bigr)^\delta$}, drawn as \tikz[baseline=-0.1cm] \draw[rho] (0.0,0) to (1,0);. (We again suppress the arbitrary orientation since it is irrelevant.) For each $P\in \mathcal{P}_N(J)$ and $\{i,j\} \in P$, it connects $v_i$ to $v_j$.
Set $I_e=1$ for these edges, as well as $(a_e,r_e)=(4+\kappa, -1)$ for some $\kappa\in (0,1)$ and call these $\rho$-type edges.
 \item Finally, there are $|J^c|$ edges, for each $i\in J^c$ one connecting $v_{k(i)}$ to $w_i$ representing the function \smash{$F^{x_{k(i)}}_\delta(x_i)$} for which we set $(a_e, r_e)=(0,0)$. (We again omit them from the schematic drawing below).
 \end{itemize}
We illustrate this with an example for $N=3$, $q=4$.
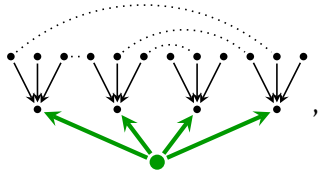
\[
\begin{tikzpicture}[scale=0.35,baseline=0.3cm]
	\node at (0,-1) [root] (root) {};
	
	\node at (-4.5,1) [dot] (a) {};
	\node at (-1.5,1) [dot] (b) {};
	\node at (1.5,1) [dot] (c) {};
	\node at (4.5,1) [dot] (d) {};	
	
	\node at (-5.5,3) [dot] (aa) {};
	\node at (-4.5,3) [dot] (ab) {};	
	\node at (-3.5,3) [dot] (ac) {};	
	
	\node at (-2.5,3) [dot] (ba) {};
	\node at (-1.5,3) [dot] (bb) {};
	\node at (-0.5,3) [dot] (bc) {};
	
	\node at (2.5,3) [dot] (cc) {};
	\node at (1.5,3) [dot] (cb) {};
	\node at (0.5,3) [dot] (ca) {};
	
	\node at (5.5,3) [dot] (dc) {};
	\node at (4.5,3) [dot] (db) {};	
	\node at (3.5,3) [dot] (da) {};	
	
	\draw[testfcn] (a) to (root);
	\draw[testfcn] (b) to (root);
	\draw[testfcn] (c) to (root);
	\draw[testfcn] (d) to (root);

	\draw[kernel] (aa) to (a);
	\draw[kernel] (ab) to (a);
	\draw[kernel] (ac) to (a);
	\draw[kernel] (ba) to (b);
	\draw[kernel] (bb) to (b);
	\draw[kernel] (bc) to (b);
	\draw[kernel] (ca) to (c);
	\draw[kernel] (cb) to (c);
	\draw[kernel] (cc) to (c);
	\draw[kernel] (da) to (d);
	\draw[kernel] (db) to (d);
	\draw[kernel] (dc) to (d);
	
	\draw[rho] (aa) to[bend left=40] (db);
	\draw[rho] (ac) to (ba);
	\draw[rho] (bb) to[bend left=40] (da);
	\draw[rho] (bc) to[bend left=40] (cb);
\end{tikzpicture}, \qquad
\]

Now we are ready to check the assumptions of the modified Hairer--Quastel criterion, Theorem~\ref{thm:slight variant of HQ}.
The first item of Assumption~\ref{ass:mainGraph_weak}, i.e., that for every edge $e\in \CCE$, one has $m_{e}-r^{-}_e<|\fraks|=4$ is clearly satisfied.
The second item, i.e., \eqref{e:assEdges}, is easily checked as well. Thus, it remains to calculate the quantities
\[
\tilde \alpha = 4|\CCV\setminus \CCV_\star| - \sum_{e\in \CCE} a_e-R(\CCG), \qquad
 R(\CCG):= \max_{\bar \CCV \subset \CCV\setminus \CCV_\star} \Biggl( 4|\bar \CCV|- \sum_{e\in \CCE(\bar\CCV) } a_e \Biggr)\vee 0.\]

Let $m$ be the number of vertices in $\CCV\setminus \CCV_\star$ which have only one incident edge. Then one sees choosing \smash{$\bar \CCV$} to be this subset maximises
$R(\CCG)= 4m - 2m= 2m$.
Since $\smash{\sum_{e\in \CCE}} a_e
= 2Nq + \smash{\frac{Nq-m}{2}(4+\kappa)}
$,
 we find that
\[
 \tilde \alpha
 = 4|\CCV\setminus \CCV_\star| - \sum_{e\in \CCE} a_e-R(\CCG)
=-(Nq-m)\frac{\kappa}{2}.
\]
Finally, the difference bound is obtained as usual.
\end{proof}

\section{Proof of Proposition~\ref{prop:varblowup}: Variance blow-up}\label{sec:proof of}
In this section, we prove Proposition~\ref{prop:varblowup}, that shows that renormalisation of the g-PAM equation with deterministic functions generically leads to mean or variance blow-up.
\begin{proof}[Proof of Proposition~\ref{prop:varblowup}]
We first note that
 in terms of the Green's function $\Gamma$ of the differential operator $\partial_t -\sum_{ij}a_{ij}(x)\partial_i\partial_j$, we have
\[\bigl(u_\delta(1,\cdot)\xi_\delta,\phi\bigr)=\int_0^1 \int_{\mathbb{R}^2}\int_{\mathbb{R}^2}\phi(x)\Gamma(s,x,y)\xi_\delta(y)\xi_\delta(x)\,\dd y \dd x \dd s.\]
As in the proof of Theorem \ref{lem:dumbbell}, let $\psi\colon C_c^\infty(B_{1/4}(0))$ taking values in [0,1] be such that $\psi= 1$ on $B_{1/8}(0)$.
For $Z$ as in \eqref{eq:frozen kernel}, recall \eqref{eq:space dep frozen}
\begin{gather*}
G^{x}(y)= \int_{0}^1 Z^{x}(t,y)\, \dd t.
\end{gather*}
Set as in
\eqref{kernels}
\begin{gather*}
K^{x}(y,z)=\psi(y-z) G^{x}(y-z), \qquad
R(x,y):= \int_{0}^1\Gamma(s,x,y)\, \dd s -K^{x}(x-y).
\end{gather*}
We then split $u_\delta(1,\cdot)$ accordingly, that is for $x\in \mathbb{R}^2$
\begin{gather*}
u_\delta(1,x)=\underbrace{\int_{\mathbb{R}^2}K^x(x,y)\xi_\delta(y)\,\dd y}_{:=\,\bar{u}_\delta(x)}+\underbrace{\int_{\mathbb{R}^2}R(x,y)\xi_\delta(y)\,\dd y}_{:=\,r_\delta(x)}.
\end{gather*}
Arguing exactly as for \eqref{EstimateErrorTermKernel}, it holds that for any $q<\infty$
\begin{gather}\label{BoundRestPartRDelta}
\sup_{\delta>0}\mathbb{E}\bigl[((r_\delta\,\xi_\delta,\phi))^q\bigr]<\infty.
\end{gather}
Expanding the square and using Young's inequality, we further have
\begin{align*}
\mathbb{E}\bigl[(u_\delta\xi_\delta-c_\delta,\phi)^2\bigr]&=\mathbb{E}\bigl[(\bar{u}_\delta\xi_\delta-c_\delta,\phi)^2\bigr]+2\mathbb{E} [(\bar{u}_\delta\xi_\delta-c_\delta,\phi)(r_\delta \xi_\delta,\phi)] + \mathbb{E}\bigl[(r_\delta \xi_\delta,\phi)^2\bigr]\\
&\stackrel{\eqref{BoundRestPartRDelta}}{\geq} \frac{1}{2}\mathbb{E}\bigl[(\bar{u}_\delta(\cdot)\xi_\delta-c_\delta,\phi)^2\bigr]+\mathcal{O}(1).
\end{align*}
Therefore, it is sufficient to show that either, there exists a $\phi$ such that
\begin{gather}\label{eq:SingMeanVarBlowUp}
\limsup_{\delta\to 0} \mathbb{E}[(\bar{u}_\delta\xi_\delta-c_\delta,\phi)]=\infty \qquad\text{or}
\qquad\limsup_{\delta\to 0} \mathbb{E}\bigl[(\bar u_\delta\xi_\delta-c_\delta,\phi)^2\bigr]=\infty.
\end{gather}

We now split the proof into two steps, showing the two alternatives in \eqref{eq:SingMeanVarBlowUp} separately.

{\it Step $1$. Mean blow-up. }We first compute $\mathbb{E}[(\bar{u}_\delta(1,\cdot)\xi_\delta,\phi)]$ using Gaussian integration by parts~\eqref{eq:GaussianIntByParts}:
\begin{align*}
&\mathbb{E}[(\bar{u}_\delta(1,\cdot)\xi_\delta,\phi)]=\int_{\mathbb{R}^2}\phi(x)\mathbb{E}[K^x(x-y)\xi_\delta(x)\xi_\delta(y)]\,\dd y \dd x\\
&\qquad=\int_{\mathbb{R}^2}\int_{\mathbb{R}^2}\phi(x)\mathbb{E}[K^x(x-y)] \bigl(\rho^{*2}\bigr)^\delta(x-y)\,\dd x \dd y\\\
&\qquad\quad +\rho^\delta* \sigma(0)\int_{\mathbb{R}^2}\int_{\mathbb{R}^2}\phi(x)\mathbb{E}\bigl[\partial_{\uparrow}^2K^x(x-y)\bigr]\rho^\delta* \sigma(x-y)\,\dd x \dd y\\\
&\qquad= \int _{\mathbb{R}^2}\phi(x)\mathbb{E}\bigl[\hat{c}_\delta^{\<Xi2>}(x)\bigr]\dd x+\rho^\delta* \sigma(0)\int_{\mathbb{R}^2}\int_{\mathbb{R}^2}\phi(x)\mathbb{E}\bigl[\partial_{\uparrow}^2K^x(x-y)\bigr]\rho^\delta* \sigma(x-y)\,\dd x \dd y,
\end{align*}
where we have set as in \eqref{kernels} for any $x\in \mathbb{R}^2$
\begin{gather}\label{DivergingTermVarianceBU}
\hat{c}_{\delta}^{\<Xi2>}(x) :=\int_{\mathbb{R}^2}K^{x}(y)\bigl(\rho^{*2}\bigr)^\delta(y)\dd y.
\end{gather}

Since the second term remains bounded, the mean blow-up occurs for any deterministic sequence $\{c_\delta\}_{\delta>0}$ which does not cancel the divergence, i.e., if there exists $\phi\in \cc\bigl(\mathbb{T}^2\bigr)$ such that as $\delta\downarrow 0$
\begin{gather}\label{FirstAmternativeMeanBU}
\int_{\mathbb{R}^2}\phi(x)c_\delta(x)\dd x\neq
\int_{\mathbb{R}^2}\phi(x)\mathbb{E}\bigl[\hat{c}^{\<Xi2>}_\delta(x)\bigr]\dd x+\mathcal{O}(1),
\end{gather}
which corresponds to the first alternative in \eqref{eq:SingMeanVarBlowUp}.

{\it Step $2$. Variance blow-up. }
We now assume that \eqref{FirstAmternativeMeanBU} is not satisfied, where without loss of generality we assume that
\smash{$
c_\delta(x)=\E\bigl[\hat{c}^{\<Xi2>}_\delta\bigr].
$}
Then, note that by the triangle inequality
\[
\E\bigl[\bigl(c_\delta-\hat{c}^{\<Xi2>}_\delta, \phi \bigr)^2\bigr]^{1/2}\leq \E\bigl[(\bar{u}_\delta\xi_\delta-c_\delta, \phi)^2\bigr]^{1/2} + \E\bigl[\bigl(\bar{u}_\delta\xi_\delta-\hat{c}^{\<Xi2>}_\delta, \phi \bigr)^2\bigr]^{1/2}.
\]
Recall that in the proof of Lemma~\ref{lem:dumbbell}, using Lemmas~\ref{lem:old_term} and~\ref{lem:new_term}, we essentially\footnote{To obtain the estimate without the additional positive renormalisation present therein, one uses the variant of the Hairer--Quastel criterion in Theorem~\ref{thm:slight variant of HQ} to obtain the appropriate analogue of Lemma~\ref{lem:old_term}.}
 already bounded the second term on the right-hand side. For the left-hand side, we use \eqref{DivergingTermVarianceBU} to get
\[\E\bigl[\bigl(\E\bigl[\hat{c}^{\<Xi2>}_\delta\bigr]-\hat{c}^{\<Xi2>}_\delta, \phi \bigr)^2\bigr]=
\int_{\mathbb{R}^2}\int_{\mathbb{R}^2}\phi(x)\phi(x') {\rm Cov}\bigl(\hat{c}_\delta^{\<Xi2>}(x),\hat{c}_\delta^{\<Xi2>}(x') \bigr)\,\dd x \dd x'. \]
Thus, using the asymptotic for any $x\in \mathbb{R}^2$
\[\hat{c}_\delta^{\<Xi2>}(x)=\int_{\mathbb{R}^2}K^x(y)\rho^\delta(y)*\rho^\delta(y)\,\dd y=-\frac{1}{2\pi}\frac{1}{\sqrt{{\rm det}(a(x))}}\log(\delta)+\mathcal{O}(1),\]
we obtain
\[
\E\bigl[\bigl(\E\bigl[\hat{c}^{\<Xi2>}_\delta\bigr]-\hat{c}^{\<Xi2>}_\delta, \phi \bigr)^2\bigr]\underset{\delta\downarrow 0}{\sim} \frac{1}{4\pi^2}\log^2(\delta)\int_{\mathbb{R}^2} \phi(x)\,\dd x\int_{\mathbb{R}^2}\phi(x')\,\dd x'{\rm Cov}\Bigl(\tfrac{1}{\sqrt{{\rm det}(a(x))}},\tfrac{1}{\sqrt{{\rm det}(a(x'))}}\Bigr).\]
To conclude, since
\[
{\rm Cov} \colon\ (x,x')\in \mathbb{R}^2\times \mathbb{R}^2\mapsto {\rm Cov}\Big(\tfrac{1}{\sqrt{{\rm det}(a(x))}},\tfrac{1}{\sqrt{{\rm det}(a(x'))}}\Big)
\]
is continuous and
\[
{\rm Cov}(0,0)=\text{Var}\Bigl(\tfrac{1}{\sqrt{{\rm det}(a(0))}}\Bigr)>0
\]
by assumption, there exists $\mathbb{A}\subset \mathbb{R}^2$ such that if $\operatorname{supp}\phi\subset \mathbb{A}$ and $\phi>0$ on the interior of $\mathbb{A}$, then
\[\int_{\mathbb{R}^2}\int_{\mathbb{R}^2}\phi(x)\phi(x'){\rm Cov}(x,x')\, \dd x'\dd x>0.
\tag*{\qed}
\]
\renewcommand{\qed}{}
\end{proof}

\appendix
\section{Gaussian integration by parts}\label{App:Gaussian int}
In this appendix, we derive Gaussian integration by parts formulae with respect to Wick products. All results of this section are `folklore' and based on Isserlis' theorem:
Given a smooth function $F \colon \mathbb{R}^n\rightarrow \mathbb{R}$ such that all derivatives have polynomial growth and real-valued centred jointly Gaussian random variables $Y$, $\{X_i\}_{i=1}^n$, it holds that
\begin{equation}\label{eq:GaussianIntByParts}
\E[F(X_1,\dots, X_n) Y ]= \sum_{i=1}^n \E[ \partial_{i}F(X_1,\dots,X_n)] \E[X_i Y].
\end{equation}
Thus, we shall always work under the following assumption.
\begin{assumption}\label{ass:FS}
Whenever working with a function $F\colon \mathbb{R}^n \to \mathbb{R}$ and random variables $\{X_\ell, Z_k\}_{k,\ell\in \mathbb{N}}$, we assume that
\begin{enumerate}\itemsep=0pt
\item For every $k\in \mathbb{N}^n$, there exist an $m\in \mathbb{N}$ such that
\smash{$\sup_{\vert x\vert>1} \frac{D^k F(x)}{|x|^m}< \infty$}.
\item The random variables $\{X_\ell, Z_k \}$ are centred and jointly Gaussian.
\end{enumerate}
\end{assumption}

Given a finite set $J$ with an even number of elements, we denote by $\mathcal{P}(J)$ the set of pairings.\footnote{Recall that a paring $P$ of $J$ is a subset $P\subset 2^{J}$ such that $|p|=2$ for all $p\in P$ and $p\cap p'= \varnothing$ whenever $p\neq p'$. We use the convention that if $J=\varnothing$ the set $\mathcal{P}(J)$ contains one element (the empty pairing).}
Then, the following lemma is a simple consequence of \eqref{eq:GaussianIntByParts}.
\begin{Lemma}
In the setting of Assumption~{\rm \ref{ass:FS}}, for $n,m\in \mathbb{N}$
\begin{align*}
&\E\Biggl[ F(X_1,\dots,X_n)\prod_{i=1}^m Z_i \Biggr] \\
&
= \!\sum_{
{\underset{\vert J\vert\ {\rm even}}
{J\subset\{1,\dots,m\}}}
}\sum_{P\in \mathcal{P}(J)}\prod_{\{i,j\}\in P}\!\mathbb{E}[Z_iZ_j]\!\sum_{k\colon J^c\rightarrow\{1,\dots,n\}}\!\mathbb{E}\biggl[\frac{\partial^{\vert J^c\vert}}{\prod_{i\in J^c}\partial_{k_i}}F(X_1,\dots,X_n)\biggr]\prod_{i\in J^c}\mathbb{E}[Z_i X_{k_i}],
\end{align*}
where the sum over $J\subset\{1,\dots,m\}$ is over even subsets and $\mathcal{P}(J)$ over pairings in $J$.
\end{Lemma}
Similarly, one easily checks the following special case of the more general Proposition~\ref{prop:wickInteByPartsGeneral}.
\begin{Lemma}\label{lem:wick_integration by parts}
In the setting of Assumption~{\rm \ref{ass:FS}}, for $n,m\in \mathbb{N}$
\begin{align*}
&\E\Biggl[ F(X_1,\dots,X_n)\prod_{i=1}^m (Z_{2i-1}Z_{2i}- \E[Z_{2i-1}Z_{2i}]) \Biggr] \\
&\qquad=
\sum_{
{\underset{\vert J\vert \ {\rm even}}
{J\subset\{1,\dots,2m\}}}
}\sum_{P\in \mathcal{P}_2(J)}\prod_{\{i,j\}\in P}\mathbb{E}[Z_iZ_j]\sum_{k\colon J^c\rightarrow \{1,\dots,n\}}\mathbb{E}\biggl[\frac{\partial^{\vert J^c\vert}}{\prod_{i\in J^c}\partial_{k(i)}}F\biggr]\prod_{i\in J^c}\mathbb{E}\bigl[Z_i X_{k(i)}\bigr],
\end{align*}
where the sum over $J$ is over even subsets of $\{1,\dots,2m\}$ and
\[\mathcal{P}_2(J):=\{ P\in \mathcal{P}(J)\mid \{2i-1,2i\}\notin P \text{ for any } i\in \{1,\dots,m\}\} \subset \mathcal{P}(J).
\]
\end{Lemma}

\subsection*{Wick products}
For any $N\geq 1$ and any jointly Gaussian (centered) random variables $\{Z_i\}_{i=1}^N$, their Wick product \smash{$\mathop{\diamond}_{i=1}^N Z_i$} is defined recursively by
\begin{gather}\label{DefWickProd}
\mathop{\diamond}_{i=1}^{N} Z_i:= Z_N \cdot \biggl( \mathop{\diamond}_{i=1}^{N-1} Z_i\biggr)- \sum_{j=1}^{N-1} \E[Z_N Z_j]
\cdot
\mathop{\diamond}_{\underset{i\neq j}{i=1}}^{N-1} Z_i,
\end{gather}
where \smash{$\mathop{\diamond}_{i=1}^{1} Z_i:= Z_1$}.
Next, we recall the following properties of the Wick product:
\begin{itemize}\itemsep=0pt
\item (Symmetry) For any permutation $\sigma\colon \{1,\dots,N\} \to \{1,\dots,N\}$, it holds that
\begin{gather}\label{eq:permutation}
\mathop{\diamond}_{i=1}^{N} Z_i = \mathop{\diamond}_{i=1}^{N} Z_{\sigma(i)}.
\end{gather}
 Therefore, this product directly extends to Gaussians indexed by any finite index set.
\item (Differential identity) For any $\ell\in \{1,\dots,N\}$, it holds
\begin{gather}\label{ResponseWickDerivative}
\frac{\partial}{\partial Z_\ell}\mathop{\diamond}_{i=1}^N Z_i=\mathop{\diamond}_{\underset{i\neq \ell}{i=1}}^N Z_i.
\end{gather}
\end{itemize}

\begin{prop}[integration by parts with respect to Wick products]\label{prop:wickInteByPartsGeneral}
In the setting of Assumption~{\rm \ref{ass:FS}}, for $n,m,N \in \mathbb{N}$
\begin{align*}
&\E\Biggl[ F(X_1,\dots,X_n)\prod_{i=0}^{m-1} \biggl(\mathop{\diamond}_{j=1}^N Z_{Ni+j}\biggr) \Biggr]\\
&=\!\!
\sum_{
{\underset{\vert J\vert \ {\rm even}}
{J\subset\{1,\dots,Nm\}}}
}\sum_{P\in \mathcal{P}_N(J)}\prod_{\{i,j\}\in P}\!\mathbb{E}[Z_iZ_j]\!\sum_{k\colon {J^c}\to \{1,\dots,n\}}\!\!\mathbb{E}\biggl[\frac{\partial^{\vert J^c\vert}}{\prod_{i\in J^c}\partial_{k_i}} F(X_1,\dots,X_n)\biggr]\prod_{i\in J^c}\mathbb{E}[Z_i X_{k_i}],
\end{align*}
where
\begin{gather*}
\begin{split}
& \mathcal{P}_N(J):=\{ P\in \mathcal{P}(J)\mid |\{Ni+1,Ni+2,\dots, Ni+N \} \cap p|\leq 1 \\
& \hphantom{\mathcal{P}_N(J):=\{} \text{for any } i\in \{0,\dots,m-1\}\text{ and any } p\in P \}.
\end{split}
\end{gather*}
\end{prop}

\begin{proof}
We show the formula by induction over $m\geq 1$ for fixed $N\geq 1$.

For the initialization $m=1$, we note that for any non-empty $J\subset \{1,\dots,N\}$, one has $\mathcal{P}_N(J)=\varnothing$. Therefore, only the term in $J=\varnothing$ contributes and the formula to prove reads
\begin{gather}
\mathbb{E}\biggl[F(X_1,\dots,X_n)\mathop{\diamond}_{i=1}^N Z_{i}\biggr]\nonumber\\
\qquad{}=\sum_{k\colon \{1,\dots,N\}\to \{1,\dots,n\}}\mathbb{E}\biggl[\frac{\partial^N}{\prod_{i\in J^c}\partial_{k_i}}F (X_1,\dots, X_n)\biggr]\prod_{i=1}^N \mathbb{E}[X_{k_i}Z_{i}].\label{Casem1IntByParts}
\end{gather}
To see this, first note that by the recursive definition of the Wick product \eqref{DefWickProd}, Isserlis' formula~\eqref{eq:GaussianIntByParts} and~\eqref{ResponseWickDerivative},
\begin{align*}
&\mathbb{E}\biggl[F(X_1,\dots,X_n)\mathop{\diamond}_{j=1}^N Z_{j}\biggr]\\
&\qquad{}\stackrel{\eqref{DefWickProd}}{=}\mathbb{E}\biggl[F(X_1,\dots,X_n)Z_N\mathop{\diamond}_{j=1}^{N-1} Z_{j}\biggr]-\sum_{\ell=1}^{N-1}\mathbb{E}[Z_{\ell}Z_N]\mathbb{E}\Biggl[F(X_1,\dots,X_n)\mathop{\diamond}_{\underset{j\neq \ell}{j=1}}^{N-1}Z_{j}\Biggr]\\
&\qquad{}\stackrel{\eqref{eq:GaussianIntByParts},\, \eqref{ResponseWickDerivative}}{=}
\sum_{\ell=1}^n\mathbb{E}\biggl[\partial_\ell F(X_1,\dots,X_n)\mathop{\diamond}_{i=1}^{N-1} Z_{i}\biggr]\mathbb{E}[X_{\ell}Z_N].
\end{align*}
Thus \eqref{Casem1IntByParts} follows by iterating.

For the induction step with $m\geq 2$, we write
\begin{gather}
\mathbb{E}\Biggl[F(X_1,\dots,X_n)\prod_{i=0}^{m}\biggl(\mathop{\diamond}_{j=1}^N Z_{Ni+j}\biggr)\Biggr]\nonumber\\
\qquad{}=\mathbb{E}\Biggl[F(X_1,\dots,X_n)\mathop{\diamond}_{j=1}^N Z_{Nm+j}\prod_{i=0}^{m-1}\biggl(\mathop{\diamond}_{j=1}^N Z_{Ni+j}\biggr)\Biggr].\label{intpartloc}
\end{gather}

Setting $X_{n+j}:= Z_{Nm+j}$ for $j\in \{1,\dots,N\}$ and writing
\[\tilde{F}(X_1,\dots,X_{n+N}):= F(X_1,\dots,X_n) \mathop{\diamond}_{j=1}^N X_{Nm+j}, \]
 by induction hypothesis \eqref{intpartloc} is equal to
\begin{gather}
\sum_{
{\underset{\vert J\vert\text{ even}}
{J\subset\{1,\dots,Nm\}}}
}\sum_{P\in \mathcal{P}_N(J)}\prod_{\{i,j\}\in P}\mathbb{E}[Z_iZ_j]\nonumber\\
\qquad{}\times\sum_{k\colon {J^c}\to \{1,\dots,n+N\}}\mathbb{E}\biggl[\frac{\partial^{\vert J^c\vert}}{\prod_{i\in J^c}\partial_{k_i}} \tilde{F}(X_1,\dots,X_{n+N})\biggr]\prod_{i\in J^c}\mathbb{E}[Z_i X_{k_i}].\label{intpartloc2}
\end{gather}
We define
\[J^{c,1}:= k^{-1}\{1,\dots,n\}, \qquad J^{c,2}:= k^{-1}\{n+1,\dots,n+N\}\]
and observe that
\[\frac{\partial^{\vert J^c\vert}}{\prod_{i\in J^c}\partial_{k_i}} \tilde{F}(X_1,\dots,X_{n+N})
=
\frac{\partial^{\vert J^{c,1}\vert}}{\prod_{i\in J^{c,1}}\partial_{k_i}} \frac{\partial^{\vert J^{c,2}\vert}}{\prod_{i\in J^{c,2}}\partial_{k_i}} \tilde{F}(X_1,\dots,X_{n+N})
\]
vanishes, except when $k|_{J^{c,2}}$ is injective in which case we write
\[{I}(J,k)= \bigl\{k(j)-n+Nm \mid j\in J^{c,2}\bigr\}, \qquad I^c(J,k)= \{Nm+1,\dots,Nm+N\}\setminus {I}(J,k)\]
and find that
\begin{align*}
&\frac{\partial^{\vert J^{c,1}\vert}}{\prod_{i\in J^{c,1}}\partial_{k_i}} \frac{\partial^{\vert J^{c,2}\vert}}{\prod_{i\in J^{c,2}}\partial_{k_i}} \tilde{F}(X_1,\dots,X_{n+N})\\
&\qquad=
\frac{\partial^{\vert J^{c,1}\vert}}{\prod_{i\in J^{c,1}}\partial_{k_i}} {F}(X_1,\dots,X_{n})
\mathop{\diamond}_{i\in \{n+1,\dots,n+N\}\setminus k(J^{c,2}) } X_{i}\\
&\qquad=
\frac{\partial^{\vert J^{c,1}\vert}}{\prod_{i\in J^{c,1}}\partial_{k_i}} {F}(X_1,\dots,X_{n})
\mathop{\diamond}_{i\in I^c({J,k})} Z_{i}.
\end{align*}
Also noting that
$\prod_{i\in J^c}\mathbb{E}[Z_i X_{k_i}]= \prod_{i\in J^{c,1}}\mathbb{E}[Z_i X_{k_i}] \prod_{i\in {J^{c,2}}} \mathbb{E}[Z_i X_{k_i}]$
this means that \eqref{intpartloc2} equals
\begin{align}\label{intbypartloc3}
&\sum_{
{\underset{\vert J\vert\text{ \tiny even}}
{J\subset\{1,\dots,Nm\}}}
}\sum_{P\in \mathcal{P}_N(J)}\prod_{\{i,j\}\in P}\mathbb{E}[Z_iZ_j]\\
&\quad\times\!
\sum_{k\colon {J^c}\to \{1,\dots,n+N\}}\!\mathbb{E}
\biggl[\frac{\partial^{\vert J^{c,1}\vert}}{\prod_{i\in J^{c,1}}\partial_{k_i}} {F}(X_1,\dots,X_{n})
\mathop{\diamond}_{i\in I^c(J,k)} Z_{i} \biggr]
\prod_{i\in J^{c,1}}\mathbb{E}[Z_i X_{k_i}] \prod_{i\in {J^{c,2}}} \mathbb{E}[Z_i X_{k_i}].\nonumber
\end{align}
Next, we use \eqref{Casem1IntByParts} to see that
\begin{align*}
&
\biggl[\frac{\partial^{\vert J^{c,1}\vert}}{\prod_{i\in J^{c,1}}\partial_{k_i}} {F}(X_1,\dots,X_{n})
\mathop{\diamond}_{i\in I^c(J,k)} Z_{i} \biggr]\\
&
\qquad{}=
\sum_{k'\colon {I^c(J,k)}\to \{1,\dots,n\} }
\mathbb{E}\biggl[
\frac{\partial^{\vert I^c(J,k) \vert}}{\prod_{i'\in I^c(J,k)}\partial_{k'_{i'}}} \frac{\partial^{\vert J^{c,1}\vert}}{\prod_{i\in J^{c,1}}\partial_{k_i}} {F}(X_1,\dots,X_{n})
\biggr]\prod_{i'\in I^c(J,k)}\mathbb{E}\bigl[Z_{i'} X_{k'_{i'}}\bigr].
\end{align*}
And thus inserting this into \eqref{intbypartloc3} yields
\begin{align}
&\sum_{
{\underset{\vert J\vert\text{ even}}
{J\subset\{1,\dots,Nm\}}}
} \sum_{P\in \mathcal{P}_N(J)}\nonumber\\
&\qquad{} \times
\sum_{k\colon {J^c}\to \{1,\dots,n+N\}} \sum_{k'\colon {I^c(J,k)}\to \{1,\dots,n\} }
\mathbb{E}\biggl[
\frac{\partial^{\vert I^c(J,k) \vert}}{\prod_{i'\in I^c(J,k)}\partial_{k'_{i'}}} \frac{\partial^{\vert J^{c,1}\vert}}{\prod_{i\in J^{c,1}}\partial_{k_i}} {F}(X_1,\dots,X_{n})
\biggr] \nonumber \\
&\qquad{}
\times \prod_{\{i,j\}\in P}\mathbb{E}[Z_iZ_j]
\prod_{i'\in I^c(J,k)}\mathbb{E}[Z_{i'} X_{k'_{i'}}]
\prod_{i\in J^{c,1}}\mathbb{E}[Z_i X_{k_i}] \prod_{i\in {J^{c,2}}} \mathbb{E}[Z_i X_{k_i}] \label{prefinal}.
\end{align}
Finally, we perform a substitution to complete the proof.
\begin{itemize}\itemsep=0pt
\item Noting that the map
$J^{c,2}\to I(J,k)$, $i\mapsto k_i -n + Nm$ is bijective,
we define $\hat{k}\colon I(J,k) \to J^{c,2}$ to be its inverse.
In particular \smash{$\prod_{i\in {J^{c,2}}} \mathbb{E}[Z_i X_{k_i}] = \prod_{i\in I(J,k)} Z_i Z_{\hat{k}_i} $}.

\item We set $\tilde{P}= P\cup \bigcup_{i\in I(J,k)} \bigl\{i, \hat{k}_i\bigr\}$ which belongs to $P_N(\{1,\dots,N(m+1)\})$.

\item We set $\tilde{J}= J \cup I(J,k) \cup \hat{k}(I(J,k))$, from which it follows that
\[
\tilde{J}^{c}(J,k):=\{1,\dots,N(m+1) \}\setminus \tilde{J}= J^{c,1} \cup I^c(J,k).
\]

\item Set $\tilde{k}\colon \tilde{J}^{c}\to \{1,\dots,n\}, \ i\mapsto \mathbf{1}_{i\leq Nm} k_i + \mathbf{1}_{i> Nm} k'_i
$
\end{itemize}
We note that for each $J$, $P$, $k$,
\begin{align}
& \mathbb{E}\biggl[
\frac{\partial^{\vert I^c(J,k) \vert}}{\prod_{i'\in I^c(J,k)}\partial_{k'_{i'}}} \frac{\partial^{\vert J^{c,1}\vert}}{\prod_{i\in J^{c,1}}\partial_{k_i}} {F}(X_1,\dots,X_{n})
\biggr]
\prod_{\{i,j\}\in P}\mathbb{E}[Z_iZ_j]
\prod_{i'\in I^c(J,k)}\mathbb{E}\bigl[Z_{i'} X_{k'_{i'}}\bigr]\nonumber\\
&\quad\qquad \times
\prod_{i\in J^{c,1}}\mathbb{E}[Z_i X_{k_i}] \prod_{i\in {J^{c,2}}} \mathbb{E}[Z_i X_{k_i}] \nonumber \\
&\qquad=
 \mathbb{E}\biggl[
\frac{\partial^{ \tilde{J}^c \vert}}{\prod_{i\in \tilde{J}^c}\partial_{\tilde{k}_i}} {F}(X_1,\dots,X_{n})
\biggr]
\prod_{\{i,j\}\in \tilde{P}}\mathbb{E}[Z_iZ_j]
\prod_{i\in \tilde{J}^c }\mathbb{E}\bigl[Z_i X_{\tilde{k}_i}\bigr].\label{finalise}
\end{align}
Finally, inserting \eqref{finalise} into \eqref{prefinal} and observing that the domains of summation agree concludes the proof.
\end{proof}

\section{The Hairer--Quastel Criterion and a Variation on it}\label{App:HairerQuastel}
Let
 $\CCG= (\CCV, \CCE)$ be a finite directed multi-graph with edges $e \in \CCE$ labelled by $(a_e, r_e)\in \mathbb{R}_+\times \{-1,0,1\}$.
 For $e\in\CCE$, we write $(e_-,e_+)$ for the pair of vertices such that $e$ is directed from the vertex $e_-$ to $e_+$.
We are given a kernel assignment $J_{e}\colon \mathbb{R}^d\times \mathbb{R}^d\setminus \triangle \to \mathbb{R}$ of compactly supported kernels.
We assume that the kernels satisfy
$\| J\|_{\alpha,\gamma}:= \sup_{|k|_\fraks <\gamma} \sup_{0<|x-y|_\fraks\leq 1} |x-y|_{\fraks}^{\alpha+|k|_{\fraks}} \bigl|D^kJ(x,y)\bigr| < +\infty $.
We then define for edges $e$ with $r_e\geq 0$
\[\hat{J}_e(x_{e_-},x_{e_+}):= J_e(x_{e_-},x_{e_+}) - \mathbf{1}_{r_e=1} J(x_{v_0},x_{e_-})\]
(where the sum is empty for $r_e=0$). Whenever $r_{e}=-1$, we assume the assignment is translation invariant and we are additionally given a real number $I_e\in \mathbb{R}$. Then, let
\[
\hat{J}_{e} (\phi) = \int J_e(x) \bigl( \phi(x)- \phi(0) \bigr)dx + I_e \phi(0),
 \]
which is well defined for any smooth compactly supported testfunctions $\phi$ whenever $a_e+r_e<|\fraks|$.

We assume $\CCG$ contains $q\geq 1$ distinguished edges\footnote{The green arrows in the schematics.} $e_{\star,1},\ldots,e_{\star,p}$ connecting $\star \in \CCV$ to distinct distinguished vertices $v_{\star,1},\ldots,v_{\star,q}$.
We write $\CCV_\star= \{ \star, v_{\star,1},\ldots,v_{\star,q}\}$ and
$\CCV_0= \CCV\setminus \{\star \}$.
We assume that whenever there are several edges connecting two vertices, at most one has non-zero renormalisation $r_e$ and that for that edge $r_e>0$.

To the distinguished edges we associate the kernel assignment $J_e=\phi^\lambda$ with label $(0,0)$.
This data provides us with a real number
\begin{gather}\label{eq:quantitiy to bound}
\hat{\mathcal{I}}_\lambda = \int_{(\mathbb{R}^d)^{\CCV_0}} \prod_{e\in \CCE} \hat{J}_e(x_{e_-}, x_{e_+}) \,\dd x.
\end{gather}

For a subset $\bar\CCV \subset \CCV$, define the following sets:
\begin{itemize}\itemsep=0pt
\item $\CCE^\uparrow (\bar \CCV) :=
 \bigl\{ e\in \CCE \mid e\cap \bar\CCV = e_- \ \text{and} \ r_e>0 \bigr\}$ (outgoing edges),
 \item $\CCE^\downarrow (\bar \CCV) = \bigl\{ e\in \CCE \mid e\cap \bar\CCV = e_+ \ \text{and} \ r_e>0\bigr\}$ incoming edges,
 \item $\CCE_0 (\bar \CCV) = \bigl\{ e\in \CCE \mid e\cap \bar\CCV = e\bigr\}$ internal edges,
 \item $\CCE (\bar \CCV) = \bigl\{ e\in \CCE \mid e\cap \bar\CCV \neq \varnothing\bigr\}$ incident edges.
\end{itemize}

\begin{assumption}\label{ass:mainGraph}
The resulting directed graph $(\CCV,\CCE)$ with labels $( a_e,r_e)$ satisfies
$r_e = 0$ for all edges connected to $\star$.
Furthermore, assume the following:
\begin{enumerate}\itemsep=0pt
\item[$1.$] For every edge $e\in \CCE$, one has $a_{e}-r^{-}_e<|\fraks|$, where $r^{-}_e:= r_e\wedge 0$.
\item[$2.$] For every subset $\bar \CCV \subset \CCV_0$ of cardinality at least $3$,
\begin{equation}\label{e:assEdges}
\sum_{e \in \CCE_0 (\bar \CCV)} a_e < \bigl(|\bar \CCV| - 1\bigr)|\fraks|.
\end{equation}
\item[$3.$] For every subset $\bar \CCV \subset \CCV$ containing $0$ of cardinality at least $2$,
\begin{equation}\label{e:assEdges2}
\sum_{e \in \CCE_0 (\bar \CCV) } a_e
+\sum_{e \in \CCE^{\uparrow} (\bar \CCV) } ( a_e+ r_e -1) - \sum_{e \in \CCE^{\downarrow} (\bar \CCV) } r_e< \bigl(|\bar \CCV| - 1\bigr)|\fraks|.
\end{equation}
\item[$4.$] For every non-empty subset $\bar \CCV \subset \CCV\setminus \CCV_\star$,
\begin{equation}\label{e:assEdges3}
 \sum_{e\in \CCE(\bar\CCV)\setminus \CCE^{\downarrow}(\bar \CCV) } a_e
 +\sum_{e\in \CCE^{\uparrow}(\bar\CCV)} r_e
- \sum_{e \in \CCE^\downarrow_+(\bar \CCV)} (r_e-1)
> |\bar \CCV||\fraks|.
\end{equation}
\end{enumerate}
\end{assumption}

\begin{Theorem}[{\cite[Theorem~A.3]{HQ18}}]\label{thm:HQ}
 Let $\CCG=(\CCV,\CCE)$ be a finite directed multigraph with labels $\{a_e,r_e\}_{e\in \CCE}$ and kernels
$\{J_e\}_{e\in \CCE}$ and real numbers $\{ I_{e} \}$
with the resulting graph satisfying Assumption~{\rm \ref{ass:mainGraph}} and its preamble.
Then, there exist $C<\infty$ depending only on the
structure of the graph $(\CCV,\CCE)$ and the value of the constants $I_e$ such that
\begin{equation*}
\hat{\mathcal{I}}_\lambda
\le C\lambda^{\tilde \alpha}\prod_{e\in \CCE} \|J_e\|_{a_e;{2}},
\end{equation*}
for $0<\lambda\le 1$, where
\begin{equation*}
\tilde \alpha = |\fraks||\CCV\setminus \CCV_\star| - \sum_{e\in \CCE} a_e.
\end{equation*}
\end{Theorem}

\subsection*{A slight variation of the criterion}
In this appendix, we derive a variant of the criterion. Given a graph and kernel assignment as above, we consider the quantity
\begin{gather*}
\mathcal{I}_\lambda = \int_{(\mathbb{R}^d)^{\CCV_0}}
\prod_{\substack{e\in \CCE\\
r_e\geq 0
}} J_e(x_{e_-}, x_{e_+})
\prod_{\substack{e\in \CCE\\
r_e<0
}} \hat{J}_e(x_{e_-}, x_{e_+}) \,\dd x.
\end{gather*}
Note that compared to \eqref{eq:quantitiy to bound} here positive edges are not `recentered'.
\begin{assumption}\label{ass:mainGraph_weak}
The resulting directed graph $(\CCV,\CCE)$ with labels $( a_e,r_e)$ satisfies
$r_e = 0$ for all edges connected to $\star$.
Furthermore, assume the following:
\begin{enumerate}\itemsep=0pt
\item[$1.$] For every edge $e\in \CCE$, one has $a_{e}-r^{-}_e<|\fraks|$.
\item[$2.$] For every subset $\bar \CCV \subset \CCV_0$ of cardinality at least $3$,
\begin{gather*}
\sum_{e \in \CCE_0 (\bar \CCV)} a_e < \bigl(|\bar \CCV| - 1\bigr)|\fraks|.
\end{gather*}
\end{enumerate}
\end{assumption}

Then one has the following variant of \cite[Theorem~A.3]{HQ18}.
\begin{Theorem}\label{thm:slight variant of HQ}
 Let $\CCG=(\CCV,\CCE)$ be a finite directed multigraph with labels $\{a_e,r_e\}_{e\in \CCE}$ and kernels
$\{J_e\}_{e\in \CCE}$ and real numbers $\{ I_{e} \}$
with the resulting graph satisfying Assumption~{\rm \ref{ass:mainGraph_weak}}.
Then, for any $\eps>0$ there exist $C<\infty$ depending only on the
structure of the graph $(\CCV,\CCE)$, the value of the constants $I_e$ and $\eps$ such that
\begin{equation*}
\mathcal{I}_\lambda
\le C\lambda^{\tilde \alpha-\eps}\prod_{e\in \CCE} \|J_e\|_{a_e;2},
\end{equation*}
for $0<\lambda\le 1$, where
\begin{equation*}
\tilde \alpha = |\fraks||\CCV\setminus \CCV_\star| - \sum_{e\in \CCE} a_e-R(\CCG), \qquad
 R(\CCG):= \max_{\bar \CCV \subset \CCV\setminus \CCV_\star} \biggl( |\bar \CCV||\fraks|- \sum_{e\in \CCE(\bar\CCV) } a_e \biggr)\vee 0.
\end{equation*}
\end{Theorem}

Since much of the proof follows ad verbatim as the one of \cite[Theorem~A.3]{HQ18}, we shall only point to the necessary adaptations in the argument therein.
A first simplification happens in~Definition~A.5 therein, where it suffices to take $\mathbf{n}\in \mathbb{N}$ since no positive renormalisation is present. Then one works with the same multiscale clustering as in \cite[Section~A.2]{HQ18}, that is consider a rooted binary tree $T$ with fixed distinguished inner node $v_*$. Denote by $T^\circ$ the set of inner nodes of $T$. Let $\mathcal{N}_{\lambda}(T^\circ)$ be all integer labellings $\ell\colon T^\circ \to \mathbb{N}$ which preserve the partial order on $T^\circ$ which has the root as minimal element, and which are such that $2^{-\ell_{v_*}}\leq \lambda$. Finally, given such a labelling $\eta\colon T^\circ \to \mathbb{R}$, we set
\[
\mathcal{I}_\lambda(\eta)= \sum_{\ell \in \mathcal{N}_{\lambda}(T^\circ)} \prod_{\nu\in T^\circ} 2^{-\ell_{\nu} \eta_{\nu}}.
\]
Next, we set
$|\eta|:= \sum_{\nu\in T^\circ} \eta_{\nu}$. Furthermore,
for any
 $\nu<\nu_*$ we write $s(\nu)$ to denote the unique minimal element of $\{ u\in T^\circ \mid \nu< u\leq\nu_* \} $ and set
 \smash{$H_\nu= \sum_{u > \nu, u\ngeq s(\nu)} \eta (u)$}. Finally, set
\begin{gather}\label{eq:HQmodeq}
R(\eta):= \max_{\nu\leq \nu_*}\biggl(\sum_{u <\nu} H_u \biggr)\vee 0.
\end{gather}
 One has the following modification of \cite[Lemma~A.10]{HQ18}, see also Remark~A.12 therein.
\begin{Lemma}\label{lem:ass_to check}
Assume that $\eta$ has the property that for every $\nu\in T^\circ$ one has $\sum_{v\geq \nu} \eta_{\nu}>0$. Then,
it holds that for any $\kappa>0$
\begin{gather*}
|\mathcal{I}_\lambda(\eta)| \lesssim \lambda^{|\eta|-R(\eta)-\kappa},
\end{gather*}
 uniformly over $\lambda\in (0, 1]$.
\end{Lemma}

With this lemma at hand one proceeds again exactly as in \cite[Section~A.3]{HQ18}, but this time setting
$ \eta(v)= |\fraks| + \sum_{e\in \hat{\CCE}} \eta_{e}(v)$ and $ \tilde{\eta}(v)= |\fraks| + \sum_{e\in \hat{\CCE}} \tilde{\eta}_{e}(v)$, where
\begin{gather}\label{eq:def eta}
 \eta_{e}(v)= -{a}_e \mathbf{1}_{e\uparrow}(v), \qquad \tilde{\eta}_{e}(v)= -{a}_e \mathbf{1}_{e\uparrow}(v)+ |r_e| \mathbf{1}_{e\in A^{-}}( \mathbf{1}_{e\uparrow}(v)- \mathbf{1}_{e\Uparrow}(v))
\end{gather}
instead of \cite[formulas (A.20) and (A.27)]{HQ18}. One checks exactly as therein that identity~(A.18) still holds since the kernels with negative regularity are translation invariant and we decomposed them as in \cite[Lemma~A.4]{HQ18}.

Thus finally, it only remains to check the following analogue of \cite[Lemma~A.19]{HQ18}.
\begin{Lemma}
The function $\tilde{\eta}$ defined in \eqref{eq:def eta} satisfies the assumptions of Lemma~{\rm \ref{lem:ass_to check}} and
$R(\eta)\leq R(\CCG)$.
\end{Lemma}

\begin{proof}
Since the first claim follows exactly as in the proof of \cite[Lemma~A.19]{HQ18}, we directly turn to the inequality.
We assume that $R(\eta)>0$, since otherwise the inequality is automatic. Let $\bar{\nu}\in T^\circ$ be the node maximising the quantity in \eqref{eq:HQmodeq},
i.e., such that
$R(\eta)=\sum_{u <\bar{\nu}} H_u $.
One then has for $U_v:=\{u\in T^\circ \mid u\ngeq v \}$
\[
\sum_{u <\bar{\nu}} H_u = \sum_{u\in U_{\bar{\nu}}}\tilde{\eta}.
\]
Denoting by $\bar \CCV$ the set of leaves attached to $U_v$, we note that $\bar \CCV \subset \CCV \setminus \CCV_*$ and
\[
\sum_{u\in U_v}\tilde{\eta}(u)= \sum_{u\in U_v}\biggl( |\fraks|+ \sum_e \tilde{\eta}_e(u) \biggr) \leq
\sum_{u\in U_v}\biggl( |\fraks|+ \sum_e {\eta}_e(u) \biggr)
=|\bar \CCV| |\fraks| -\sum_{e\in \hat{\CCE}(\bar \CCV)} {a}_e.
\tag*{\qed}
\]
\renewcommand{\qed}{}
\end{proof}

\section{Kernel spaces}\label{App:Kernels}
Here we recall the kernel spaces which quantify \cite[Assumption~5.1]{Hai14} at finite regularity.
We first introduce notations for Taylor expansions.
For $\gamma>0$ and $f\colon \mathbb{R}^{d}\to \mathbb{R}$ compactly supported such that $D^{k}f$ exists whenever $|k|_\fraks<\gamma$, let
\[
\opP_{z_{0}}^\gamma [f](z):= \sum_{|k|_\fraks<\gamma} \frac{D^{k}f(z_0)}{k!} (z-z_0)^k.
\]
Similarly, for $\gamma'>0$ and a compactly supported function of two variable $F\colon \mathbb{R}^{d}\times \mathbb{R}^d \to \mathbb{R}$ such that
$D^{k_1}D^{k_2}F$ exists whenever $|k_1|_\fraks<\gamma$, $|k_2|_\fraks<\gamma'$,
 let
\[\opP^{(\gamma,\gamma')}_{({z}_0, \bar{z}_0)}[F](z, \bar{z}) := \sum_{|k|<\gamma, |l|<\gamma' } \frac{D_1^{k} D_2^{l} F({z}_0, \bar{z}_0) }{k!\ l!} (z-z_0)^k(\bar{z}-\bar{z}_0)^l.\]
The following is \cite[Defenition~1]{Sin25}, we refer to said work for further discussion of its content.
\begin{Definition}
Let $K\colon \mathbb{R}^{d+1}\times \mathbb{R}^{d+1}\setminus \triangle \to \mathbb{R}$ be a kernel
which can be decomposed as $K(z,z')= \sum_{n\geq 0} K_n(z,z')$ where each $K_n$ is supported on
\smash{$ \bigl\{(z,z')\mid \|z-z'\|_\fraks \leq 2^{-n+1} \bigr\}$}.
On the space of such kernels $\pmb{\mathcal{K}}$, we define for $L, R \in \mathbb{R}_+$,
the norm $\|K\|_{\beta;L,R}$ as the smallest number $C$ such that there exists a decomposition $K= \sum K_n$ for which the following bounds are satisfied:
\begin{itemize}\itemsep=0pt
\item For any multi-indices $k_1$, $k_2$ satisfying $|k_1|_\fraks < L$, $|k_2|_\fraks< R$ and $n \in \mathbb{N}$, it holds that
\begin{equation*}
 \sup_{z,z'}\bigl|D_1^{k_1}D_2^{k_2} K_n(z,z') \bigr| \leq C2^{(|\fraks|-\beta +|k_1|_\fraks +|k_2|_\fraks )n}
\end{equation*}
as well as
\begin{align*}
& \sup_{z,z',\bar{z}'} \frac{\bigl|D^{k_1}_1K(z, z')- \opP^{R}_{\bar{z}' }\bigl[D^{k_1}_1K(z,\cdot ) \bigr](\bar{z}')\bigr| }{|z'-\bar{z}'|^R } \lesssim C2^{n(|\fraks|-\beta +|k_1|_\fraks+ R )},\\
& \sup_{z,z',\bar{z}} \frac{\bigl|D^{k_2}_1K(z, z')- \opP^{R}_{\bar{z}' }\bigl[D^{k_2}_2K(\cdot,z' ) \bigr](\bar{z})\bigr| }{|z-\bar{z}|^L } \lesssim C2^{n(|\fraks|-\beta +|k_2|_\fraks+ L )}
\end{align*}
and
\begin{gather*}
 \sup_{z,\bar{z},z',\bar{z}'} \frac{\bigl|K_n(z,z') -\opP^L_{\bar{z}} [K_n(\,\cdot\,, {z}')](z) -\opP^{R}_{\bar{z}'}[K_n({z}, \,\cdot\, )](z')+ \opP^{(L,R)}_{(\bar{z}, \bar{z}')}[K_n](z, {z}')\bigr| }{|z-\bar{z}|_\fraks^\gamma |z'-\bar{z}'|_\fraks^{\gamma'} }\\
 \qquad\leq C2^{n(|\fraks|-\beta +L+ R )}
\end{gather*}
uniformly over $n\geq 0$.

\item Let $\mathbf{k}_{<R}:= \bigl\{k \in \mathbb{N}^d \mid |k|_{\fraks}<R \bigr\}$ and denote by $\partial\mathbf{k}_{<R}:= \{
k\notin \mathbf{k}_{<R}\mid k-e_{\min\{ j\mid k_j\neq 0 \}\in \mathbf{k}_{<R} }
\}
$ its boundary in the sense of \cite[Appendix~A]{Hai14}. Then, for any $k_1\in \mathbf{k}_{<R}$ and $k_2\in \mathbf{k}_{<R}\cup\partial\mathbf{k}_{<R}$
\begin{equation*}
\biggl|\int_{\mathbb{R}^d} (z-z')^{k_1} D^{k_2}_2 K_n(z,z') \,\dd z\biggr| \leq C2^{-\beta n }
\end{equation*}
uniformly over $n\geq 0$ and $z'\in \mathbb{R}^{d+1}$.
\end{itemize}
Finally, we define the vector space \smash{$\bfK^\beta_{L,R}= \{ K\in \bfK \mid \|K\|_{\beta;L,R}<+\infty\}$} equipped with the norm ${\|\cdot\|_{\beta;L,R}}$.
\end{Definition}

\subsection*{Acknowledgements}
Both authors would like to thank the Oberwolfach Research Fellows (OWRF) program for supporting a research stay during which a part of this work was carried out, as well as Rhys Steele and Lucas Broux for valuable discussions during that stay. HS gratefully acknowledges financial support from the Swiss National Science
Foundation (SNSF), grant number 225606. Finally, we are grateful to the anonymous referees for their detailed reading of the article and the thereby resulting improvements.

\pdfbookmark[1]{References}{ref}
\LastPageEnding

\end{document}